\documentclass[11pt]{amsart}

\usepackage{epsf}
\usepackage{geometry}
\usepackage{epsfig}
\usepackage{amssymb}
\usepackage{amsmath}
\usepackage{amsthm}
\usepackage{graphicx}
\usepackage{relsize}
\usepackage{xcolor}
\usepackage{cases}
\usepackage{subfigure}

\usepackage{lineno}
\usepackage{marginnote}
\usepackage{rotating}
\usepackage{color}
\usepackage{bm}

\usepackage{amsfonts}
\usepackage{latexsym}
\usepackage{bm}
\usepackage{array}
\usepackage{mathrsfs}
\usepackage{flafter}
\usepackage{tcolorbox}
\definecolor{darkgreen}{rgb}{0., 0.5, 0.}
\definecolor{wite}{rgb}{1., 1., 1.}
\definecolor{darkblue}{rgb}{0., 0., 0.38}
\definecolor{brown}{rgb}{0.6, 0.3, 0.0}
\definecolor{darkred}{rgb}{0.8, 0.0, 0.0}
\definecolor{lightblue}{rgb}{0.0, 0.5, 0.9}
\definecolor{softgreen}{rgb}{0.0, 0.8, 0.}
\definecolor{applegreen}{rgb}{0.55, 0.71, 0.0}
\definecolor{brickred}{rgb}{0.8, 0.25, 0.33}
\definecolor{brightmaroon}{rgb}{0.76, 0.13, 0.28}
\definecolor{brightgreen}{rgb}{0.4, 1.0, 0.0}
\definecolor{brightpink}{rgb}{1.0, 0.0, 0.5}
\definecolor{brightturquoise}{rgb}{0.03, 0.91, 0.87}
\definecolor{brightube}{rgb}{0.82, 0.62, 0.91}
\definecolor{forshade}{rgb}{0.98, 0.9, 0.93}

\newcommand{\R}[1]{{\rm I\!R}^{#1}}
\newcommand\be{\begin{eqnarray*}}
\newcommand\ee{\end{eqnarray*}}
\newcommand\ben{\begin{eqnarray}}
\newcommand\een{\end{eqnarray}}
\def\NNN{{\boldsymbol{|}\!\!\boldsymbol{|}\,}}
\def\om{\mathlarger{\mathlarger{\omega}}}
\def\Rd{{\mathbb R}^d}
\def\tanh{{\rm tanh}}
\def\bp{\boldsymbol{p}}
\def\bq{\boldsymbol{q}}
\def\br{\boldsymbol{r}}
\def\by{\boldsymbol{y}}
\def\bn{\boldsymbol{n}}

\def\dvg{{\rm div}}

\newtheorem{ass}{assumption}

\newtheorem{assumption}[ass]{A}
\newtheorem{remark}{Remark}[section]
\newtheorem{algorithm}{Algorithm}[section]
\newtheorem{theorem}{Theorem}[section]
\usepackage{tikz,pgfplots}
\usepackage{tkz-graph}
\usetikzlibrary{shapes.arrows,decorations.markings}
\usetikzlibrary{automata}
\usetikzlibrary{positioning}

\begin{document}
% Do not eliminate the following size-changing command, please:
\large
\title{A posteriori error estimates for domain decomposition methods}
\author{J. Kraus}
\address{Faculty of Mathematics, University of Duisburg-Essen, 45127 Essen, Germany}
\email{johannes.kraus@uni-due.de}
\author{S. Repin}
\address{St. Petersburg Department of V.A. Steklov institute of Mathematics of Russian Academy of Sciences, 191023, Fontanka 27, St. Petersburg, Russia}
\email{repin@pdmi.ras.ru}
\maketitle
\begin{abstract}
Nowadays, a posteriori error control methods have formed a new important part of the numerical analysis.
Their purpose is to obtain computable error estimates in various norms and error indicators
that show distributions of global and local errors of a particular numerical solution.

In this paper,
%we briefly discuss the history of error estimation methods in the context of models based on
%partial differential equations (PDEs).
we focus on a particular class of domain decomposition methods (DDM), which are among the most efficient
numerical methods for solving PDEs.
We adapt functional type a posteriori error estimates and construct a special form of error majorant which
allows efficient error control of approximations computed via these DDM by performing only subdomain-wise 
computations. The presented guaranteed error bounds use an extended set of admissible fluxes which arise 
naturally in DDM.
\end{abstract}

\tableofcontents

\section{Introduction} \label{sec:introduction}

%%%%%%%%%%%%%%%%%%%%%%%%%%%%%%
\subsection{Domain decomposition methods}

%\subsubsection{The Schwarz alternating method and its development} 
The iteration method of Schwarz \cite{Schwartz} suggested for analysis of conformal mappings
has generated a branch of computational methods, which are nowadays among the most used
in practice for the numerical solution of PDEs. In application to elliptic boundary value problems,
proofs of convergence were first presented in the paper by S. Mikhlin~\cite{Mikhlin} and in the
book of Kantorovich and Krilov~\cite{KK}, implementation issues first discussed in~\cite{Stoutemyer1973}. 

The pioneering works on the Schwarz alternating method by
Pierre-Louis Lions~\cite{Lions_1978,Lions_1987,Lions_1989}, Matsokin and
Nepomnyaschikh~\cite{MatsokinNepomnyashchikh1981,MatsokinNepomnyashchikh1985}
provided the basis of the abstract Schwarz theory, for multiplicative methods see
also~\cite{Bramble_etal1991}, and herewith layed the foundations of modern domain
decomposition methods (DDM). Schwarz methods exist as overlapping and non-overlapping
methods, where representatives of either category can be formulated as additive or
multiplicative subspace correction, as well as hybrid methods, see, e.g.~\cite{dd_book}.

%\subsubsection{Domain decomposition, multigrid, and hybrid methods}
Broad over\-views and newer developments in the area of DDM were presented
in~\cite{dd_book,M08,DoleanJolivetNataf2015}. Together with multigrid/algebraic multigrid (MG/AMG)
methods, see, e.g.,~\cite{Ha03,TOS01,V08}, DDM have become one of the most successful classes
of iterative solution methods when it comes to cost- and energy-efficient computations of approximate
solutions of PDEs. Their description and convergence analysis in the abstract framework of subspace
correction (SSC) methods has been presented by Xu and Zikatanov~\cite{XZ02}.

There exist also hybrid approaches to construct solvers and preconditioners for a broad range of discrete
models based on PDE. The idea to combine domain decomposition and multigrid techniques can be found
as early as in~\cite{K89}. The auxiliary space multigrid (ASMG) method presented in~\cite{KLM15}
implements this concept by means of auxiliary space preconditioning, see~\cite{Xu96}, and achieves
herewith robustness with respect to general coefficient variations~\cite{KLLMZ16}. 

%\subsubsection{Substructuring methods}
Nonoverlapping DDM or substructuring methods, as they are sometimes referred to, come as primal and
dual methods. Balancing domain decomposition (BDD), see~\cite{Ma93}, is a primal method operating on
the common interface degrees of freedom (DOF) whereas finite element tearing and interconnecting (FETI)
enforces equality of DOF on substructure interfaces by Lagrange multipliers and thus is a dual method.

Nonoverlapping DDM, including BDDC (BDD based on constraints), see~\cite{Do03}, and FETI-DP
(FETI dual-primal), see~\cite{FLLPR01}, can be formulated and analyzed in a common algebraic framework,
see~\cite{MD03,MDT05,MS07}. The BDDC method enforces continuity across substructure interfaces by a
certain averaging operator.

%\subsubsection{Robust coarse spaces}
A key tool in the analysis of overlapping DD methods is the Poincar\'{e} inequality or its weighted analog
as for problems with highly varying coefficients. It is well-known that the weighted Poincar\'{e} inequality
holds only under certain conditions, e.g., in case of quasi-monotonic coefficients, see~\cite{Sa94}.
More recently the robustness of DD methods has also been achieved for problems with general coefficient
variations using coarse spaces that are constructed by solving local generalized eigenvalue problems,
see, e.g.,~\cite{EGLW12,GE10a,SVZ11}. 

The additional constraints can be viewed as subspace corrections involving the coarse basis functions,
which are subject to energy minimization. From this point of view, the BDDC method~\cite{Do03} has
also a high degree of similarity with the ASMG method~\cite{KLM15}. However, contrary to BDDC,
the latter, in general exploits overlapping subdomains where coarse DOF are associated not only with
subdomain interfaces (or boundaries) but also with their interior. The recursive application of the two-level
method then results in a full multilevel/multigrid algorithm. Its main characteristics are, contrary to standard
(variational) multigrid algorithms, an auxiliary-space correction step instead of classical coarse-grid
correction, and a coarse-grid operator arising from an additive approximation of the Schur complement
of the original system, see also~\cite{Kr06,Kr12}.

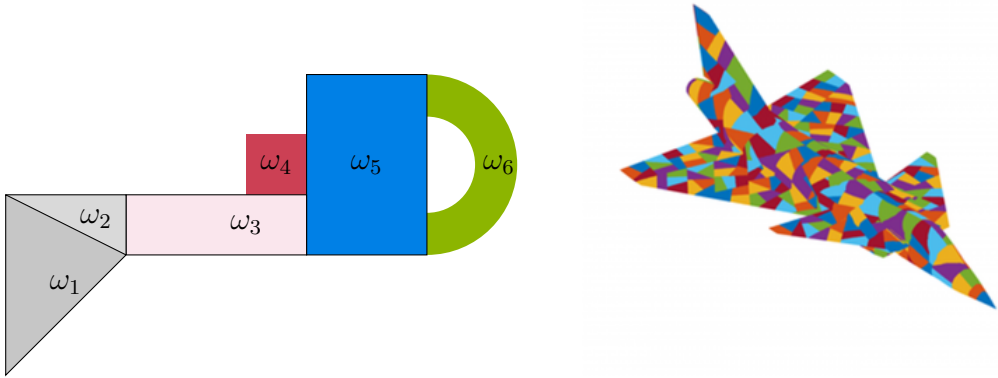
\begin{figure}[h]
\centering 
\begin{tikzpicture}[scale=0.8]
%\coordinate (4) at (1,1);
\coordinate (b)  at (1,1); % node can be labeled by a letter (a,b,...,p,q,r,s,t)
% draw a circle filled with color
\fill [brickred] (0,8) rectangle (1,9);  %  color is given in [...]
\fill [applegreen] (3,8.5) circle(1.5);
\fill [white] (3,8.5) circle(0.8);
\draw[fill=lightblue] (3,10) rectangle  (1,7);
\draw[fill=forshade] (-2,8) rectangle  (1,7);

\draw[fill=gray!30]    (-2,7) -- ++(0,1) -- ++(-1,0) -- ++(-1,0);
%\draw[fill=gray!45]    (-2,1) -- ++(-2,+1) -- ++(0,-2) -- ++(0,2);

\draw  [fill=gray!45]  (-2,7) -- ++(-2,1) -- ++(0,-3) -- ++(2,2);
%\draw[fill=white] (-1.5,1) rectangle  (-1.5,1);
\node at (-3,6.5){$\omega_1$};
\node at (-2.5,7.6){$\omega_2$};
\node at (0,7.5){$\omega_3$};
\node at (0.5,8.5){$\omega_4$};
\node at (2,8.5){$\omega_5$};
\node at (4.2,8.5){$\omega_6$};
\end{tikzpicture}
\;
\hspace{2ex}\includegraphics[width=65mm]{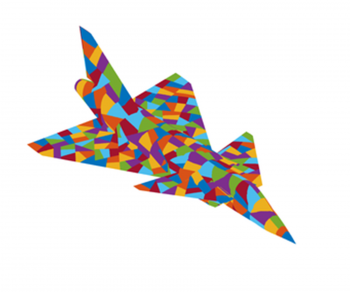}
\vspace{4ex}
\caption{Decomposition of two-- and three--dimensional domains.}
\label{figure0}
\end{figure}

%\subsubsection{Reducing complexity}
The main advantage of domain decomposition methods is their {\em ability 
to solve geometrically complicated problems by solving problems on simple
subdomains}. These methods have proven their high efficiency both for relatively
simple cases where the number of subdomains is small (as in the left picture of
Figure~\ref{figure0}) and for engineering problems related to complicated 3D
domains (right picture of Figure~\ref{figure0}).

A posteriori error estimates are not only a key element of reliable computer simulations but
are also well-suited to be combined and interact with modern SSC resulting in even more cost-
and energy-efficient computational methods for PDEs. 
\vspace{1ex}

Our work is inspired by the idea to supply DDM with guaranteed error control following 
the guiding principle,
namely, {\em reducing the estimation of global errors to error estimates for subdomains}.
\vspace{0.5ex}

%%%%%%%%%%%%%%%%%%%%%%%%%%%%%%
\subsection{Fully guaranteed a posteriori error estimates}

Difficulties related to generation of fully computable, guaranteed, and consistent
error bounds for partial differential equations are easy to explain with the paradigm
of the simplest elliptic problem $\Delta u+f=0$ in a bounded Lipschitz domain $\Omega$
with the boundary condition $u=0$ on $\partial\Omega$. The solution $u$ is an element
of the Sobolev space $H^1_0(\Omega)$ (that contains square integrable functions,
which vanish on the boundary and have square integrable
generalized derivatives of the first order). It is formally defined by the integral
identity
\ben
\label{id1}
\int\limits_\Omega(\nabla u\cdot \nabla w-fw)dx=0\qquad \forall w\in H^1_0(\Omega).
\een
Let $v\in H^1_0(\Omega)$ be an approximation of $u$
(e.g., computed by  a Ritz-Galerkin finite element method). From (\ref{id1}), it follows
that
\ben
\label{id2}
\int\limits_\Omega\nabla (u-v)\cdot \nabla w dx=\int\limits_\Omega(\nabla v\cdot\nabla w-fw)dx.
\een
The right hand side of (\ref{id2}) can be estimated via the residual of the equation
$R(v)=\Delta v+f$ in two ways. The first way is possible if $v$ has an extra regularity,
so that $R(v)\in L^2(\Omega)$. Then we integrate by parts, use the Friedrichs inequality
with a positive constant $C_F$ (which depends on $\Omega$) and arrive at the estimate
\ben
\label{id3}
\|\Delta (u-v)\|_\Omega\leq C_F\|\Delta v+f\|_\Omega,
\een
where $\|\cdot\|_\Omega$ stands for the norm of scalar and vector valued functions in $L^2(\Omega)$.
Notice that the quantity in the left hand side of (\ref{id2}) can be considered as a measure
of the error $e=u-v$. This quantity majorates $\|\nabla e\|$ with some constant multiplier (which for a
convex domain $\Omega$ is equal to 1). 

The right hand side of (\ref{id3}) is directly computable, so that at first
glance the estimate  looks attractive. However, from a computational point 
of view, this estimate has a fatal drawback: it is not applicable to a sequence of approximations
$\{v_k\}$ converging to $u$ in the energy space $H^1_0(\Omega)$. Even if the sequence is regularized
so that $R(v_k)\in L^2(\Omega)$ and the norms in (\ref{id3}) are computable and finite for any $k$,
we are unable to guarantee that the strong norm of the residual tends to zero. Hence, (\ref{id3})
does not possess an important {\em consistency} property.

Another way of error estimation operates with a weak norm of $R(v)$
\ben
\label{id4}
\sup\limits_{w\in H^1_0(\Omega)}
\frac{\int\limits_\Omega\nabla (u-v)\cdot \nabla w dx}{\|\nabla w\|_\Omega}=
\sup\limits_{w\in H^1_0(\Omega)}
\frac{\int\limits_\Omega(\nabla v\cdot\nabla w-fw)dx}
{\|\nabla w\|_\Omega}=:\NNN R(v)\NNN_{-1,\Omega}.
\een
It is easy to show that the left hand side of (\ref{id4}) coincides
with $\|\nabla e\|_\Omega$. Moreover, the norm $\NNN R(v_k)\NNN_{-1,\Omega}$
tends to zero if $v_k$ tends to $u$ in $H^1_0(\Omega)$ and, therefore, this error
relation is consistent. But here we are faced with another difficulty.
Unlike integral type norms, the norm $\NNN R(v)\NNN_{-1,\Omega}$
 is incomputable because the supremum
is taken over an infinite amount of test functions. If the set of functions
is reduced to some finite dimensional subspace, then the upper bound
in (\ref{id4}) may be lost.
An attempt to overcome
this difficulty and get a computable majorant of $\NNN R(v)\NNN_{-1,\Omega}$
using special
properties of $v$ (Galerkin orthogonality) is known in the literature as the
explicit residual method (e.g., see \cite{Ve3}). In practice, 
this way leads to
the above discussed  
class of error
indicators.

A posteriori error estimates, which are guaranteed, do not contain mesh--dependent constants,
and are valid for any conforming approximation error, are derived by more sophisticated methods
(see a consequent exposition in \cite{Repin2000,ReGruyter,MaNeRe}). Estimates of this type are
derived by purely functional methods without attraction of an information about the method and
mesh used to compute a numerical solution. Therefore, they are known as {\em a posteriori estimates
of the functional type}. These estimates do not use special properties of approximations (e.g., Galerkin
orthogonality or superconvergence) as well as of the exact solution (additional regularity). The simplest
estimate of this class for the problem \eqref{eq:1a}--\eqref{eq:1b} reads
\begin{equation}
\label{eq:4}
\Vert \nabla(u-v) \Vert_A \le \Vert A \nabla v-\by\Vert_{A^{-1}} +C_F \Vert \dvg \by+f\Vert_{\Omega},
\end{equation}
where   $C_F$ is the constant in the Friedrichs inequality for the doman $\Omega$
and we use the notation
\be
\|\bq\|^2_A:=\int\limits_\Omega A\bq\cdot\bq dx\quad\text{and}\quad
\|\bq\|^2_{A^{-1}}:=\int\limits_\Omega A^{-1}\bq\cdot\bq dx\quad\forall \bq\in [ L^2(\Omega)]^d.
\ee
 This
estimate is valid for any function $v\in H^1(\Omega)$ that satisfies the boundary condition
\eqref{eq:1c} and any vector valued function $\by$ in the space
$$
H(\Omega,\dvg):=\left\{\bq\in[ L^2(\Omega)]^d: \, \dvg \bq\in L^2(\Omega) \right\}.
$$
The function $\by$ can be viewed as an approximation of the exact flux $\bp=A\nabla u$.
It is easy to see that if $\by=\bp$ then the estimate coincides with the exact error.
If $\bq\in Q_f$, see (\ref{Qf}), then we arrive at the well known hypercircle
estimate
(see \cite{Sy,Mi})
\begin{equation}\label{eq:3}
\Vert \nabla(u-v)\Vert_A \le \inf_{\bq\in Q_f}\Vert A\nabla v -\bq\Vert_{A^{-1}}.
\end{equation}
There is a whole class of a posteriori methods based on the use of~(\ref{eq:3})
and its analogs for different boundary value problems. In them, $\by$ is
defined by post--processing of the numerical flux $A\nabla v$ (which does not belong to $Q_f$)
 and the main efforts are focused on satisfying (\ref{eq:3})
  (exactly or approximately). Here, we refer to
\cite{Ke,LaLe} and subsequent papers of many authors devoted to a posteriori methods
 based on equilibration of fluxes (e.g., see \cite{BrSh}).
 Numerical procedures used to construct a suitable equilibrated flux
can be rather
complicated, and, most importantly, they do not always guarantee that condition 
 (\ref{eq:3}) is exactly
satisfied. In the latter case, the estimate does not provide a guaranteed
 upper bound of the error, but usually serves as a good error indicator.
It should be outlined that $Q_f$ is a rather narrow subset of $H(\Omega,\dvg)$.
The space $H(\Omega,\dvg)$ admits
simple conforming approximations that preserve the continuity of $\by\cdot n$ 
on interior boundaries
(e.g., the so-called Raviart-Thomas elements). Therefore,
a numerical flux $A\nabla v$  can be  projected to $H(\Omega,\dvg)$ by
very simple post--processing procedures. Moreover, it is 
useful to make several iterations of minimization with respect
to $\by$. After that  the right hand side
of (\ref{eq:4}) gives a good and fully guaranteed majorant of the error.
In this paper, we show that this approach can be adapted to domain decomposition
methods.

% Interest is driven by its impact on natural and life sciences, economics, engineering where many 
% processes are described in terms of PDE.
% Discretization techniques such as finite difference, finite volume, and finite element methods
% reduce the continuous problem to a discrete problem which finally is expressed in form of one or
% more systems of linear algebraic equations.

In conclusion of this overview, a few words should be said about reliable control 
of the accuracy of numerical solutions to integral equations.
They arise in a number of important mathematical models, and in addition,
there are methods for analyzing 
differential equations by reducing them to integral equations (e.g., Picard--Lindel\"of
method \cite{Nev}).
Fully guaranteed a posteriori error estimates for this class of problems are based
upon the estimates derived by A. Ostrovskii for abstract iteration procedures~\cite{Os}.
The reader can find a systematic discussion of these questions in Chapter~6 of~\cite{MaNeRe}.
\vspace{1ex}

%%%%%%%%%%%%%%%%%%%%%%%%%%%%%%
\subsection{Outline of the paper and main results}

The present paper presents a new methodological approach that allows to control global error resulting in the
process of applying a DD method to solve iteratively an elliptic boundary-value problem,
e.g.~Problem~\eqref{eq:1a}--\eqref{eq:1b}, using only computations on subdomains, that is, utilizing operators which are defined
and act locally, on subdomains, only. Consequently, any entirety of approximate solutions of subdomain problems
yields explicit estimates of local (subdomain) errors as well as of the global error at any stage of the iterative
process and thus gives a picture of the error contributions associated with individual subdomains. 

The remainder of this paper is organized as follows. Section~\ref{sec:basic_problem} contains the formulation
of a simple elliptic model problem, introduces some notation, and recalls the basic structure of functional a
posteriori error estimates for an elliptic model problem. Section~\ref{sec:ddm} provides algorithmic details and
summarizes some known facts about domain decomposition methods as they can be combined with local
(subdomain-wise) a posteriori error estimation that will be analyzed hereafter in Section~\ref{sec:error_bounds}.
The main theroretical result of this work is the a posteriori error estimate stated in Theorem~\ref{thm:main_result},
which is especially designed for domain decomposition methods. Finally, Section~\ref{sec:numerics} discusses
numerical results, which, on the one hand, support the theoretical estimates presented in Section~\ref{sec:ddm},
and, on the other hand, suggest strategies and possibilities to further refine the proposed approach.

%%%%%%%%%%%%%%%%%%%%%%%%%%%

\section{Domain decomposition methods for the basic elliptic problem}\label{sec:basic_problem}

In the next section, we discuss guaranteed and fully computable error estimates
adapted to DDM using the following
 elliptic boundary-value problem
 \begin{subequations}\label{eq:1}
\begin{eqnarray}
\label{eq:1a}
&&
\dvg  \bp+\,f  = 0~~~ \text{in}~~ \Omega ,  \\
\label{eq:1c}
&& \bp=A \nabla u,\\
&&u  = u_g\;\; \text{on}~~ \Gamma =\partial{\Omega} , 
\label{eq:1b} 
\end{eqnarray}
\end{subequations}
where $\Omega \in \R{d}$ is a polygonal domain, 
$u_g$ is a given function in $H^1(\Omega)$, and (\ref{eq:1b}) is understood as
equality of the corresponding traces. 

We assume that $A$ is a symmetric positive definite (SPD) $d{\times}d$ matrix, satisfying the estimate
\begin{equation}\label{pd_A}
C_{\min} \Vert \xi \Vert^2 \le A\xi \cdot \xi \le C_{\max} \Vert \xi \Vert^2, \quad \forall \xi \in \R{d}.
\end{equation}
The corresponding weak formulation of the problem reads:
find $u \in H^1_0(\Omega)$ such that
\begin{equation}\label{eq:2}
a(u,w):=\int_{\Omega} A\nabla u\cdot \nabla w \,dx =
\int_{\Omega} f w\, dx, \  \forall w \in V_0:=H^1_0(\Omega).
%= \langle F,w\rangle \quad \mbox{for all } w \in H^1_0(\Omega)
\end{equation}
Hence the exact flux $\bp=A\nabla u \in Q_f$, where
\ben
\label{Qf}
Q_f:=\left\{ \bq\in[ L^2(\Omega)]^d: \int_{\Omega}\bq\cdot \nabla w = \int_{\Omega}fw \, dx \quad 
\forall \ w\in V_0 \right\}.
\een
\label{sec:ddm}

The iterative solvers we study here fall into the category of alternating Schwarz type methods based on an overlapping
domain decomposition.
In order to describe the basic setting of such methods formally we use the following assumptions and notation.

Let $\Omega$ be partitioned into a collection of "basic" (relatively simple, convex)
subdomains $\om_k$ (e.g., triangles or quadrilaterals for $d=2$ and tetrahedra or hexahedra for $d=3$) so that
\begin{equation}\label{eq:partition_basic}
\overline{\Omega}=\bigcup_{k=1}^{N} \overline{\om}_k\qquad
{\rm and}\quad {\om}_i \cap {\om}_j = \emptyset\;{\rm if}\;i \ne j.
\end{equation}
In addition to $\{\om_k\}$, we consider another set of Lipschitz subdomains
$\Omega_j$, $j=1,2,\ldots,M$ that are utilized by a domain
decomposition method. For the decomposition of $\Omega$ into 
$\{ \Omega_j: j=1,2,\ldots,M \}$ we  assume that 
\begin{equation}\label{eq:decomposition1}
\overline{\Omega}=\bigcup_{j=1}^{M} \overline{\Omega}_j,
\end{equation}
and
\begin{equation}\label{eq:decomposition2}
\overline{\Omega_k \cap \Omega_l}=\begin{cases}\emptyset, \\ \text{certain collection of } \overline{\om}_i, \end{cases}
\end{equation}
which implies that each subdomain $\Omega_j$ itself is partitioned into a certain set of basic subdomains, i.e.,
\begin{equation}\label{eq:partition_subdomain}
\overline{\Omega}_k = \bigcup_{j\in \mathcal{I}_k} \overline{\om}_j,
\end{equation}
where $\mathcal{I}_k$ denotes the index set corresponding to the nonoverlapping partitioning of the subdomain
$\Omega_k$ into basic subdomains $\om_j$. Note that this setting allows to consider both overlapping and
nonovelapping domain decomposition methods.

\begin{figure}[h]
\centering 
\subfigure[Two overlapping subdomains $\Omega_1$ (vertically elongated rectangle)
and $\Omega_2$ (horizontally elongated rectangle)]{\label{fig:a}
\includegraphics[width=65mm]{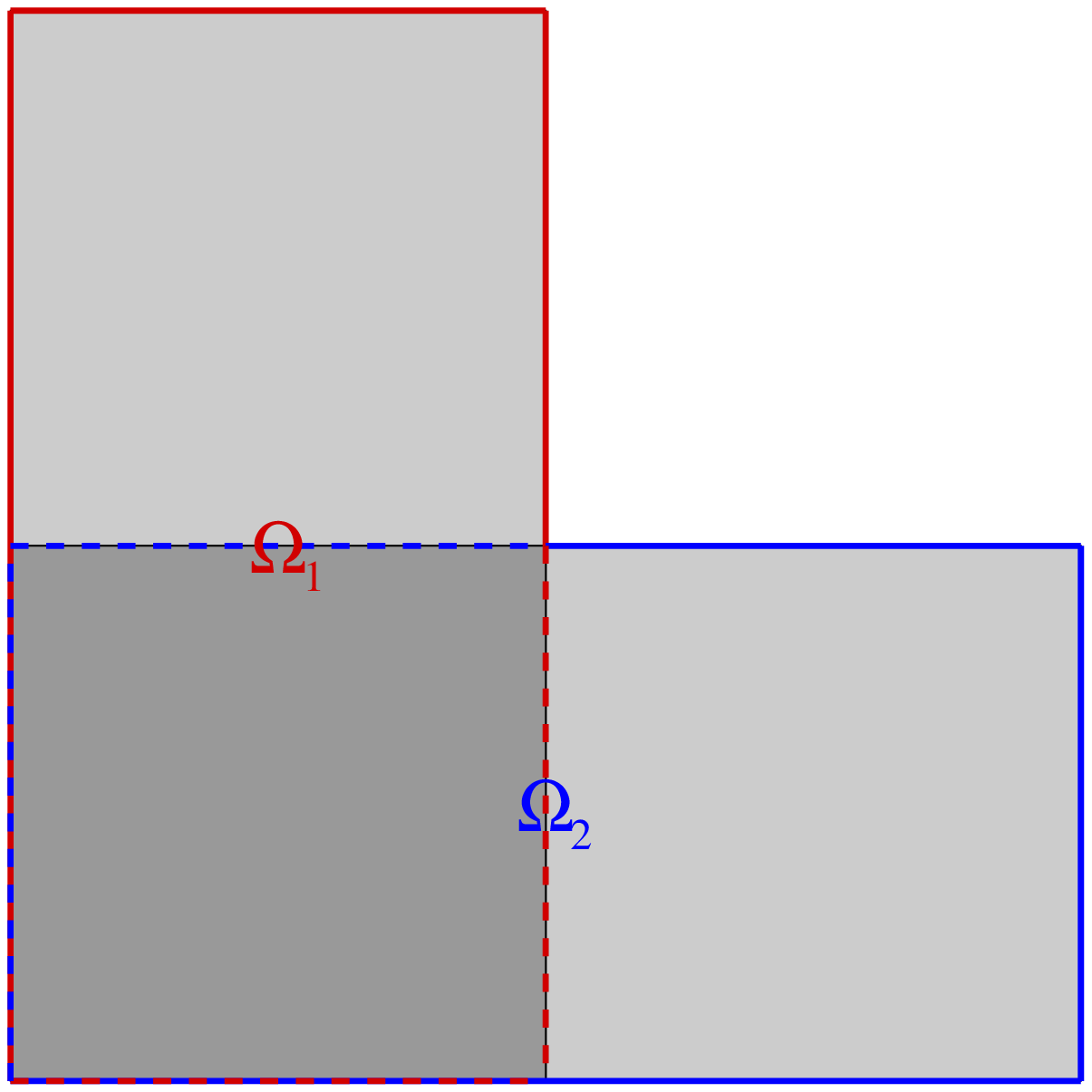}}\hspace{6ex}
\subfigure[Three non-overlapping basic subdomains $\om_1$, $\om_2$, and $\om_3$]{\label{fig:b}
\includegraphics[width=65mm]{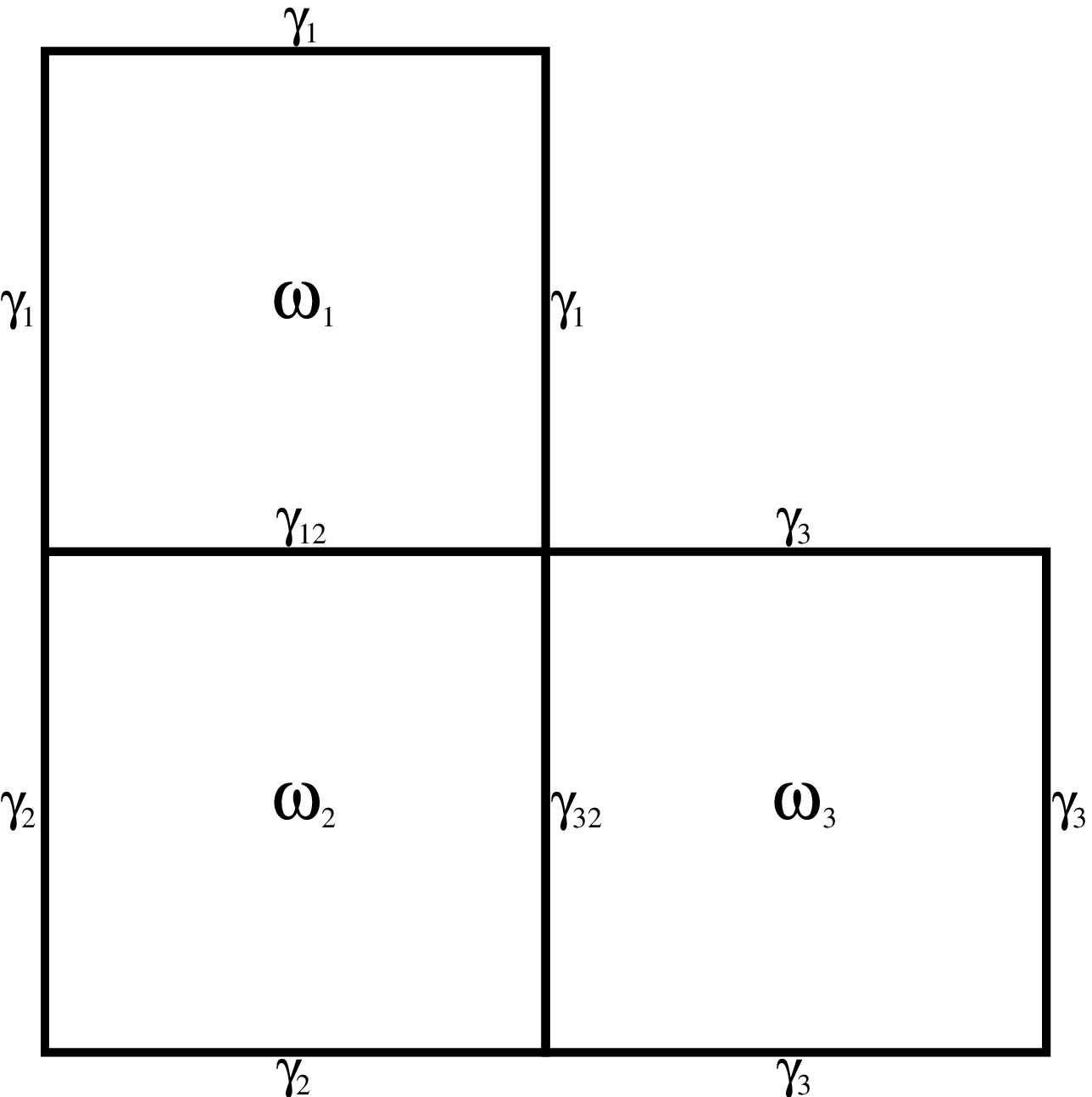}}\\
\caption{Notation for subdomains and boundaries}\label{figure1}
\end{figure}

First, let $\gamma_{ji}=\gamma_{ij}$ denote the interface of 
$\overline{\om}_i$ and $\overline{\om}_j$, that is, 
$$
\gamma_{ij} := \overline{\om}_i \cap \overline{\om}_j .
$$
We assume that $\gamma_{ij}$ has a positive
surface measure, otherwise ${\om}_i$ and ${\om}_j$
are considered as non--intersecting cells.
Next, the boundary faces (edges) $\gamma_i$ of $\overline{\om}_i$ are denoted by 
$$
\gamma_i : = \overline{\om}_i \cap \Gamma.
$$
Like this the boundary $\partial \Omega_k$ is given by 
$$
\partial \Omega_k :=
\left(
\bigcup_{i \in \mathcal{I}_k \land j \notin \mathcal{I}_k
\lor
i \notin \mathcal{I}_k \land j \in \mathcal{I}_k}
\gamma_{ij} \right) \cup
\left(
\bigcup_{j \in \mathcal{I}_k} \gamma_{j} \right)
 =:  
 \left( \bigcup_{i \notin \mathcal{I}_k} \Gamma_{k,i} \right) \cup \Gamma_k
$$ 
and $\Gamma_k=\partial \Omega_k \cap \Gamma$.
We will consider here the ``overlapping case'', that is, the case in which $\gamma_{ij} \subset \partial \Omega_k$
for $i \in \mathcal{I}_k$ and $j \notin \mathcal{I}_k$ or $i \notin \mathcal{I}_k$ and $j \in \mathcal{I}_k$  implies that
there exists $\Omega_\ell$ such that $i \in \mathcal{I}_\ell$ and $j \in \mathcal{I}_\ell$. Note that we do not care
about the order of the indices $i$ and $j$ in $\gamma_{ij}$, that is, $\gamma_{ij}$ is identified with $\gamma_{ji}$
so that for ``$i \in \mathcal{I}_k \land j \notin \mathcal{I}_k \lor i \notin \mathcal{I}_k \land j \in \mathcal{I}_k$''
we could equivalently write ``$i \notin \mathcal{I}_k \land j \in \mathcal{I}_k$'', for example.
Figure~\ref{figure1} illustrates the notation we are using on a very simple example. The left picture shows 
two overlapping subdomains $\Omega_1$ and $\Omega_2$, which are composed of two basic subdomains each, i.e.,
$\overline{\Omega}_1=\overline{\om}_1 \cup \overline{\om}_2$, $\overline{\Omega}_2=\overline{\om}_2 \cup \overline{\om}_3$.
Their intersection is given by $\Omega_1 \cap \Omega_2 = \om_2$. The right picture illustrates the three
non-overlapping basic subdomains $\om_1$, $\om_2$, and $\om_3$.\footnote{In case of a finite element discretization each basic subdomain
$\om_j$ will in general consist of a number $N_j$ of elements $T_{ji}$ of the mesh partition $\mathcal{T}_h := \{ T_{ji}: j=1,2,\ldots,N; i=1,2,\ldots,N_j \}$.}
The boundaries of $\Omega_1$ and $\Omega_2$ in this case are given by
$\partial \Omega_1 = \Gamma_{1,3} \cup \Gamma_1$ with $\Gamma_{1,3}=\gamma_{32}$ and $\Gamma_1=\gamma_1 \cup \gamma_2$
and
$\partial \Omega_2 = \Gamma_{2,1} \cup \Gamma_2$ with $\Gamma_{2,1}=\gamma_{12}$ and $\Gamma_2=\gamma_2 \cup \gamma_3$.

There are various ways to extend the classical Schwarz alternating method from two to more than two (overlapping) subdomains.
We follow reference~\cite{Lions_1987} and denote by $V_k$ the closed subspace of $H_0^1(\Omega)$ consisting of all elements
of $H_0^1(\Omega_k)$ extended by $0$ to $\Omega$. 

Then we choose some initial guess $u_0^1 \in H_0^1(\Omega)$ for the solution $u$ of problem~
\eqref{eq:1a}--\eqref{eq:1b}.
 Without loss of
generality, we may assume that $g=0$, which we can always achieve by subtracting from the solution of the inhomogeneous
equation a known function $u_g\in H^1(\Omega)$ satisfying the boundary condition $u_g=g$ on $\Gamma$ and modifying the
right hand side of equation~\eqref{eq:1a} accordingly. 

We use the following domain decomposition algorithm for solving the problem~\eqref{eq:1}, see~\cite{Lions_1987}:
\begin{algorithm}\label{alg1}
\begin{tabbing}
\hspace{2cm}\= \quad \= \quad \= \quad \= \quad \= \\[1ex]

\> \verb!Choose initial guess!  $u_0^1 \in H_0^1(\Omega)$ \\
\> \verb!for n=1,2,...! \\
\> \> \verb!for k=1 to M! \\
\> \> \> \verb!Solve! \\[-6ex] 
\end{tabbing}
\begin{subequations}\label{eq:dd}
\begin{eqnarray}
- \dvg A \nabla u_k^n &=& f~~~~ \text{in} \; \; \; \Omega_k , \hspace{2cm} \label{eq:dda} \\
u _k^n &=& u_{k-1}^n \; \; \text{on }~~~ \partial{\Omega}_k , \label{eq:ddb} 
\end{eqnarray}
\begin{tabbing}
\hspace{2cm}\= \quad \= \quad \= \quad \= \quad \= \\[-3ex]
\> \> \> \verb!Extend! \\[-5ex] 
\end{tabbing}
\begin{equation}
u _k^n \ \ \verb!to! \ \ \Omega \ \ \verb!by! \ \ u_{k-1}^{n} \label{eq:ddc}
\end{equation}
%\end{array}
\end{subequations}
\end{algorithm}

The above algorithm is the classical Schwarz alternating method, which is a multiplicative Schwarz (successive subspace correction) method
with an error propagation operator (see \cite{dd_book})
\begin{equation}\label{eq:err_prop_AS}
E_{\rm{mu}}=(I-P_M)(I-P_{M-1})\ldots (I-P_1),
\end{equation}
where $P_k=R_k^T \widetilde{P}_k$ and $\widetilde{P}_k:V\rightarrow V_k$ are defined by 
$$
a_k(\widetilde{P}_k u,v_k)=a(u, R_k^T v_k)\quad \text{for all }v_k\in V_k,\, k=1,\ldots,M. 
$$
Here $V\subset H_0^1(\Omega)$ is a finite dimensional approximation space, e.g., conforming finite element space, that is used
to discretize the global problem~\eqref{eq:2} whereas $V_k$ denote corresponding spaces of functions with local support $\bar{\Omega}_k$
in which the solutions $u_k$ of the subdomain problems 
\begin{equation}\label{weak_form}
\begin{array}{rcl}
a_k(u_k,v_k)&:=&\int_{\Omega_k}(A\nabla u_k)\cdot \nabla v_k\, dx \\
&=&\int_{\Omega_k}fv_k \,dx-a_k(\bar{u}_{g_k},v_k)=:\langle F_k,v_k \rangle \quad \forall v_k\in V_k
\end{array}
\end{equation}
are sought for and $\bar{u}_{g_k} := u_{k-1}^n \in H^1(\Omega_k)$ are  functions satisfying the subdomain boundary conditions~\eqref{eq:ddb}.
The interpolation operators $R_k^T:V_k \rightarrow V$ provide the decomposition 
\begin{equation}\label{decomposition}
V=\sum_{k=1}^M R_k^T V_k
\end{equation}
of the space $V$.

Note that it is also possible to initialize the method with an initial guess $u_{k}^{0}$ on each subdomain and then impose the
boundary conditions $u _k^n = u_{k}^{n-1}$ on $\partial{\Omega}_k$ in step~\eqref{eq:ddb} for all $k$ in which case the method
becomes an additive Schwarz (parallel subspace correction) method; step~\eqref{eq:ddc} in this case can be removed.

The convergence analysis of multiplicative and additive Schwarz methods is covered
by the abstract Schwarz theory that goes back to Pierre-Louis Lions~\cite{Lions_1978,Lions_1987,Lions_1989}, see also~\cite{dd_book}. 
 
Applying the Schwarz method according to Algorithm~\ref{alg1} we are solving local problems in $\Omega_1, \Omega_2, \ldots, \Omega_M$.
The respective numerical solutions obtained in the $n$th step are denoted by $v_1^n, v_2^n, \ldots, v_M^n$ with
$v_1^0, v_2^0, \ldots, v_M^0$ being some initial guesses.
 
We assume that for any $n=1,2,\ldots$, it holds for all $k=1,\ldots, M$
\begin{equation}
v_k^n=v_{\ell}^{n-1} \text{ on } \Gamma_{k,i} \text{ for all } i\notin \mathcal{I}_k: \exists \ell \text{ with }
\Gamma_{k,i}\subset \Omega_{\ell}. \label{eq:13a} 
\end{equation}
Hence, after the $n$th step we can introduce the conforming global approximation
\begin{equation}\label{glob_approx}
\overline{v}^n:=v_j^{n} \;\;\;\text{ in }\overline{\om}_j. 
\end{equation}

Our aim is to control the accuracy of $v_1^n, \ldots, v_M^n$ by using local majorants. This leads to respective calculations
in the subdomains $\Omega_1, \ldots ,\Omega_M$ from which we obtain respective fluxes $\by_1^n, \ldots, \by_M^n$ with
which we want to guarantee that the local problem has been solved sufficiently accurately. 

The ultimate goal is to deduce guaranteed bounds of the global error based on the ``local'' functions $v_k^n,\by_k^n$ and
on numerical computations performed on the subdomains $\Omega_k$ exclusively, $k=1, \ldots, M$.

It is well known (\cite{Lions_1987,Lions_1989,dd_book}) that
\begin{equation}\label{eq:conv_on_subdomain}
u_k^n \ \overset{n \rightarrow \infty}{\longrightarrow} \ u\vert_{\Omega_k} \mbox{ in } H^1(\Omega_k)
 \mbox{ for all } k=1,2, \ldots,M,
\end{equation}
where $u_k^n$ are the exact solutions of the subdomain problems~\eqref{eq:dd}.
To be more precise, the following convergence result holds true:
\begin{theorem}\cite[Theorem I.2]{Lions_1987}\label{thm:Lions}
The sequence $(u_m)_{m\ge1}$ defined by $u_{M(n-1)+k}:=u_k^n$, $n \ge 1$, $k=1,2,\ldots,M$, generated
by Algorithm~\ref{alg1} converges to $u$ in $V$.
In addition, if \, $\sum_{k=1}^M V_k = V$, there exists $\rho \in [0,1)$ such that 
$$
\Vert u_m-u \Vert_V \le \rho^m \Vert u_0-u \Vert_V, \quad \mbox{for all } m \ge 0.
$$
\end{theorem}

The subdomain problems~\eqref{eq:dd} on $\Omega_1, \ldots, \Omega_M$ are assumed to be much simpler to solve, which in
numerical computations is usulally due to the fact that the number of degrees of freedom (DOF) to approximate $u_k^n$ by $v_k^n$
is much smaller than the number of DOF used in the process of approximating $u$ by $v$. Typically it is also possible to choose the
subdomains $\Omega_k$ in such a way that their shape is much simpler than that of $\Omega$, which in some special cases might
even allow to use analytical methods for (approximately) solving the problems~\eqref{eq:dd}.

This means that we want to reduce a problem that is hard to solve, i.e., requires a huge number of DOF for its numerical solution with
a certain desired accuracy and/or is posed in a complicated domain, to a sequence of discrete problems of much smaller dimension in
simple domains $\Omega_k$ where we can apply very efficient solvers. Hence, in the $n$th step of the algorithm, we compute for
$k=1,2,\ldots,M$
$$
\begin{array}{ll}
v_k^n \ldots & \text{an approximation of the $n$th exact subdomain solution } u_k^n \text{ in } \Omega_k, \\[1ex]
\by_k^n \ldots & \text{a vector-function approximating the exact local flux } \bp_k^n = A \nabla u_k^n \text{ in } \Omega_k.
\end{array}
$$
The $M$ pairs $(v_k^n, \by_k^n)$, $k=1,2, \ldots ,M$, approximating respective solutions $(u_k^n,\bp_k^n)$, which we would
have on the step $n$ if local subproblems would be solved exactly, are indeed known and can be used in a posteriori error control. 

%%%%%%%%%%%%%%%%%%%%%%%%%
\section{Guaranteed bounds of errors}\label{sec:error_bounds}

To control the accuracy of the approximations obtained by DDM, one can always use the global
estimates (\ref{eq:4}) and (\ref{eq:3}). However, in this way some technical difficulties arise.
Assume that the last iteration was focused on getting a new approximation
in the subdomain $\Omega_k$. In $\Omega\setminus\Omega_k$ the approximate solution $v$ and flux $\by$
have not been changed and the corresponding errors remain the same. Hence, an efficient error control
procedure should consider only the part associated with $\Omega_k$. But in this case we cannot guarantee
continuity of the normal component of the flux across $\partial\Omega_k$ (which is required in~(\ref{eq:4})). 
To keep this continuity,
we need some global averaging procedure, which is not limited to $\Omega_k$, but also changes fluxes
in the neighbouring domains. Below we deduce a special form of the error majorant, which minimises
difficulties of this kind.

%%%%%%%%%%%%%%%%%%%%%%%%%
\subsection{A posteriori error estimate adapted to  DDM}\label{subsec:esimate}

\subsubsection{Main assumptions}
Our goal now is to deduce fully guaranteed error estimates for approximations generated by the DD method, which satisfies the assumptions~\eqref{eq:partition_basic}--\eqref{eq:partition_subdomain}.
In addition, we assume that on any step $m$ of the iteration process, the corresponding approximation satisfies the boundary condition
 $u_m=u_g$. Certainly we can use the estimate \eqref{eq:4} directly.
 However, this global estimate does not account the specifics of DDM
 approximations that change the function in one subdomain only.
 Hence it is reasonable to have an estimate in such a form that
 allows us to recompute only the part related to the lastly considered
 subdomain and utilise all other data computed on previous iterations.
 The key problem arising in this concept is related
 to proper regularity (conformity) of the approximations
 obtained by DDM. With the above mentioned conditions,
 any approximation $u_m$ will be $H^1$--conforming, so that the 
 conformity problem
 is related to approximations of the flux, which may not belong to the
 space $H(\Omega,\dvg)$. This may happen because the continuity of normal components
 on $\gamma_{kj}:=\overline\om_k\cap\overline\om_j$ is not guaranteed. Hence we introduce
 a ``broken $H(\dvg)$ space''
 \begin{equation*}
H(\bigcup_j \om_j,\dvg):=\left\{ \bq \in[ L^2(\Omega)]^d: \bq=\bq_j\;{\rm in}\;\om_j,\;\bq_j\in H(\om_j,\dvg)\;
 \mbox{ for all } \om_j \subset \Omega \right\}.
\end{equation*}
The space $H(\bigcup_j \om_j,\dvg)$ contains vector valued functions
that have square integrable divergence only locally, in subdomains
$\om_j$. It is much larger than $H(\Omega,\dvg)$ and in a sense
too large to be used as the set of possible fluxes $\by$ because these
vector valued functions should satisfy at least some weak continuity of
normal components on $\gamma_{kj}$.
Therefore, we  restrict the set of admissible
fluxes that are further used in a posteriori estimates.
\begin{assumption}
\label{ass1}
Let  $\by \in H(\bigcup_j \om_j,\dvg)$ satisfy the condition
\begin{equation}\label{eq:D4}
\{ \by_k \cdot \bn_{kj}-\by_j\cdot \bn_{kj} \}_{\gamma_{kj}}=0 \quad \forall \gamma_{kj}
\end{equation}
where $\by_k=\by\vert_{\overline{\om}_k}$ and $\by_j=\by\vert_{\overline{\om}_j}$, $\gamma_{kj}$ is a common boundary of $\om_k$ and $\om_j$
and $\bn_{kj}$ is normal to $\gamma_{kj}$. Here 
$
\{ r \}_{\gamma_{kj}} := \frac{1}{\vert \gamma_{kj} \vert} \int_{\gamma_{kj}} r \, ds.
$
\end{assumption}
This condition means that the mean value of the ``normal flux jump'' is zero on the interface $\gamma_{kj}$, i.e., the continuity of normal
flux is satisfied in a weak (integral) sense. 

Also, we impose one more condition.
\begin{assumption}\label{ass2}
Let the local fluxes be weakly equilibrated, that is,
\begin{equation}\label{eq:D6}
\{ \dvg \by_j +f\}_{\om_j}=0
\end{equation}
where
$
\{ r \}_{\om_j} := \frac{1}{\vert \om_j \vert} \int_{\om_j} r \, dx.
$
\end{assumption}
The conditions \eqref{eq:D4} and \eqref{eq:D6} are easy to satisfy by a suitable correction
of the numerical flux (a question which is addressed later).

%%%%%%%%%%%%%%%%%%%%%%%%%%%
\subsubsection{Derivation of  error majorants}
First, we rewrite the integral identity~\eqref{eq:2} in the form
\begin{multline}
\int_{\Omega} A\nabla (u-v) \cdot \nabla w \, dx = \int_{\Omega} \left(fw-A\nabla v\cdot \nabla w \right) dx\\
\label{eq:D2}
=\sum_{j=1}^N \int_{\om_j} \left(fw-A\nabla v\cdot \nabla w \right) dx .
\end{multline}
We use the divergence theorem
\begin{equation}\label{eq:D5}
\int_{\om_j} (\dvg \by_j w+ \by_j \cdot \nabla w) \, dx =
\int_{\om_j} \dvg (\by_j w) \, dx=\int_{\partial \om_j} \by_j\cdot \bn w \, ds
\end{equation}
 and arrive at the relation
\begin{equation}\label{eq:basic_rep}
\begin{array}{cl}
\displaystyle\int_{\Omega} A\nabla (u-v) \cdot \nabla w \, dx & =
\displaystyle\sum_{j=1}^N\displaystyle \int_{\om_j}\left[(\dvg \by_j+f)w+(\by_j-A\nabla v)\cdot \nabla w \right] \, dx 
\\ & -\displaystyle\sum_{j=1}^N \displaystyle\int_{\partial \om_j}\by_j \cdot \bn w \, ds .
\end{array}
\end{equation}
In view of  \eqref{eq:D6}, we estimate the first summand of the integral
over $\om_j$
as follows:
\begin{equation}\label{eq:local_est_1}
 \displaystyle \int_{\om_j} (\dvg \by_j+f)w \, dx = \displaystyle\int_{\om_j}(\dvg \by_j+f)(w-c_j) \, dx
\le \Vert \dvg \by_j +f\Vert_{\om_j} \Vert w-c_j\Vert_{\om_j},
\end{equation}
where $c_j=\{w\}_{\om_j}$. Then by the Poincar\'{e} inequality
\begin{equation}\label{eq:local_est_2}
\Vert w - c_j \Vert_{\om_j} \le C_{{\rm P},j} \Vert \nabla w\Vert_{\om_j}
\end{equation}
with local constants $C_{{\rm P},j} = C_{{\rm P},j} (\om_j)$ we obtain
\begin{equation}\label{eq:global_est_1}
\sum_{j=1}^N \int_{\om_j}(\dvg \by_j+f) w \, dx \le
\left(\sum_{j=1}^N  C_{{\rm P},\max}^2 \Vert \dvg \by_j+f \Vert^2_{\om_j} \right)^{1/2} \Vert \nabla w\Vert_{\Omega},
\end{equation}
where
\ben
\label{poincare}
C_{{\rm P},\max} := \max\limits_{j=1,2,...,N} C_{{\rm P},j}
\een
For the second summand of the integral we use the inequality
\begin{equation}\label{eq:D7}
\int_{\om_j} (\by_j-A\nabla v)\cdot \nabla w \, dx \le \|\by_j-A\nabla v\|_{A^{-1},\om_j} \|\nabla w\|_{A,\om_j} ,
\end{equation}
where $\|\by\|^2_{A,\om_j}:=\int\limits_{\om_j}A\by\cdot\by\,dx$.
By summation over all subdomains $\om_j$ and applying the discrete Cauchy-Schwarz inequality, we obtain
\begin{eqnarray}
\displaystyle\sum_{j=1}^N \int_{\om_j} (\by_j-A\nabla v)\cdot \nabla w \,dx  &\le&
\left(\displaystyle\sum_{j=1}^N \Vert \by_j-A\nabla v\Vert_{A^{-1},\om_j}^2\right)^{1/2} \left(\displaystyle\sum_{j=1}^N \Vert \nabla w\Vert_{A,\om_j}^2\right)^{1/2} \nonumber \\
&=& \left(\displaystyle\sum_{j=1}^N \Vert \by_j-A\nabla v\Vert_{A^{-1},\om_j}^2\right)^{1/2} \Vert \nabla w\Vert_{A,\Omega} . \label{eq:global_est_2}
\end{eqnarray}
Now we turn to the last sum in~\eqref{eq:basic_rep} which expands as
\begin{equation}\label{eq:last_sum}
\sum_{j=1}^N \int_{\partial \om_j} \by_j \cdot \bn \, w \, ds = \sum_{\gamma_{kj} \in\mathcal{E}_1}\int_{\gamma_{kj}} (\by_j\cdot \bn_{kj}-\by_j\cdot\bn_{kj}) \, w \, ds +
\sum_{\gamma_j \in\mathcal{E}_2}\int_{\gamma_j} \by_j \cdot \bn \, w \, ds 
\end{equation}
where $\mathcal{E}_1$ denotes the set of all interior interfaces and $\mathcal{E}_2$ denotes the set of all boundary faces, i.e., $\gamma_j \in\mathcal{E}_2$ if and only if
$\gamma_j \cap \Gamma = \gamma_k$.
For the case of
full Dirichlet boundary conditions \eqref{eq:1b}, the test function $w$
is equal to zero on $\Gamma$
and, therefore,
$$
\int_{\partial \om_k\cap \Gamma} \by_k \cdot \bn \, w \, ds=0 .
$$
Hence the last integral in the right hand
side of \eqref{eq:global_est_2} vanishes. It is worth noting,
that in the case of mixed boundary conditions $\Gamma=\Gamma_D\cup\Gamma_N$, where
$\Gamma_N \ne \emptyset$, the Assumption~\ref{ass1} must be appended with the condition 
\ben
\label{eq:D6plus}
\{\by_k\cdot \bn\}_{\gamma_k \cap \Gamma_N}=0\; \text{for all}\; \gamma_k
\;\text{such that}\;
\gamma_k \cap \Gamma_N \ne \emptyset.
\een
Then the last integral also vanishes and the error majorant can be derived by the same arguments as presented below.

To get an upper bound for~\eqref{eq:last_sum} we have to estimate the quantity
\begin{equation}\label{eq:global_est_3}
 \sum_{\gamma_{kj} \in\mathcal{E}_1}\int_{\gamma_{kj}} (\by_k\cdot \bn_{kj}-\by_j\cdot\bn_{kj}) \, w \, ds .
\end{equation}
For that purpose we need special (Poincar\'{e} type) inequalities of the form
\begin{equation}\label{eq:D8}
\begin{array}{rcl}
\Vert w\Vert_{\gamma} &:=&
\Vert w\Vert_{L^2(\gamma)} \\
&\le& C_{\rm P}(\gamma,\omega) \Vert \nabla w \Vert_{\omega}, \quad \forall w\in \tilde{H}^1(\gamma,\omega):=\{v\in H^1(\omega): \{v\}_{\gamma}=0\},
\end{array}
\end{equation}
where $\omega$ is a bounded Lipschtz domain in $\Rd$ and $\gamma$
is a connected part of the boundary $\partial\omega$ having positive surface measure.
Sharp constants for the inequality~\eqref{eq:D8} (and for some other analogous inequalities) 
have been derived in~\cite{Nazarov_2015}.
In particular, if $d=2$, $\omega=(0,h_1)\times(0,h_2)$, and $\gamma=\{x_1=0,x_2=[0,h_2]\}$, then
$$
C_{\rm P}(\gamma,\omega)=\left(\frac{\pi}{h_2}
\tanh \left(\frac{\pi}{h_2}\right)\right)^{-1/2}.
$$

\begin{remark} An important property of \eqref{eq:D8} is the monotonicity with respect to expansion of the domain. It is easy to see that if $\omega_1\subset \omega_2$ 
and for both domains $\gamma\subset \partial \Omega_i$, $i=1,2$ then the constant $C_{\rm P}(\gamma,\omega_1)$ can be used in the estimate for $\omega_2$ as well, i.e., 
$C_{\rm P}(\gamma,\omega_2)\le C_{\rm P}(\gamma,\omega_1)$. 
For this reason we can use sharp constants derived for some basic (relatively simple) domains as upper bounds of the constants related
to more complicated domains.
\end{remark}

Using A1 and the H\"older's inequality we obtain
\begin{equation}\label{eq:local_est_3}
\begin{array}{rl}
 \displaystyle\int_{\gamma_{kj}}(\by_k \cdot \bn_{kj}-\by_j \cdot \bn_{kj})w\, ds & = \displaystyle\int_{\gamma_{kj}}(\by_k \cdot \bn_{kj}-\by_j \cdot \bn_{kj})(w-c)\, ds \\[3ex]
 & \le \Vert \by_k \cdot \bn_{kj} - \by_j \cdot \bn_{kj}\Vert_{\gamma_{kj}} \Vert w-c\Vert_{\gamma_{kj}}.
\end{array}
\end{equation}
Now we use~\eqref{eq:D8} to estimate the last norm in~\eqref{eq:local_est_3}. Notice that $\gamma_{kl}$ belongs
to two neighbouring subdomains. Therefore, we can estimate the
norm in two ways:
$$
\Vert w-c\Vert_{\gamma_{kj}}\le C_{\rm P}(\gamma_{k,j},\om_j)\Vert \nabla w \Vert_{\om_j}
$$
and
$$
\Vert w-c\Vert_{\gamma_{kj}}\le C_{\rm P}(\gamma_{k,j},\om_k)\Vert \nabla w \Vert_{\om_k}.
$$
By summation we get
$$
2\Vert w-c\Vert^2_{\gamma_{kj}}\le\left[ C_{\rm P}^2(\gamma_{k,j},\om_j)+C_{\rm P}^2(\gamma_{k,j},\om_k) \right] \Vert \nabla w \Vert^2_{\om_j\cup\,\om_k}
$$
and hence
$$
\Vert w-c\Vert_{\gamma_{kj}}\le \underbrace{\sqrt{\frac{C_{\rm P}^2(\gamma_{k,j},\om_j)^2+C_{\rm P}^2(\gamma_{k,j},\om_k)^2}{2}}}_{=:\beta_{kj}} \Vert \nabla w \Vert_{\om_j\cup\,\om_k}.
$$
Returning to~\eqref{eq:local_est_3}, we obtain
$$
\left \vert \int_{\gamma_{kj}} (\by_k \cdot \bn_{kj}-\by_j \cdot \bn_{kj}) w \, ds\right\vert \le \Vert (\by_k-\by_j)\cdot \bn_{kj} \Vert_{\gamma_{kj}}
\beta_{kj}\Vert \nabla w\Vert_{\om_k \cup \om_j}.
$$
Summation over all $\gamma_{kj}\in \mathcal{E}_1$ yields the estimate
$$
\begin{array}{l}
\displaystyle\sum_{\gamma_{kj}\in \mathcal{E}_1} \Vert (\by_k-\by_j)\cdot \bn_{kj}\Vert_{\gamma_{kj}} \beta_{kj} \Vert \nabla w \Vert_{\om_j\cup\om_k} \\
\le \left(\displaystyle\sum_{\gamma_{kj}\in\mathcal{E}_1} \Vert (\by_k-\by_j)\cdot \bn_{kj}\Vert_{\gamma_{kj}}^2\beta_{kj}^2 \right)^{1/2} 
\left( \displaystyle\sum_{\gamma_{kj}\in\mathcal{E}_1}\Vert \nabla w\Vert^2_{\om_k\cup\om_j} \right)^{1/2}.
\end{array}
$$
Now we denote by
$E_{\max}$ the maximum number of interfaces associated with a single basic subdomain $\om_k$.
Then we conclude from the above estimate that 
\begin{equation}\label{eq:global_est_4}
\begin{array}{l}
\displaystyle\sum_{\gamma_{kj}\in \mathcal{E}_1} \Vert (\by_k-\by_j)\cdot \bn_{kj}\Vert_{\gamma_{kj}} \beta_{kj} \Vert \nabla w \Vert_{\om_j\cup\om_k} \\
\le \sqrt{E_{\max}} 
\left(\displaystyle\sum_{\gamma_{kj}\in\mathcal{E}_1} \Vert (\by_k-\by_j)\cdot \bn_{kj}\Vert_{\gamma_{kj}}^2\beta_{kj}^2 \right)^{1/2} \Vert \nabla w\Vert_{\Omega}.
\end{array}
\end{equation}
Collecting the estimates~\eqref{eq:global_est_2},~\eqref{eq:global_est_1} and~\eqref{eq:global_est_4} yields 
$$
\begin{array}{rl}
 \displaystyle\int_{\Omega} A\nabla(u-v)\cdot \nabla w\, dx &
 \le \left( \displaystyle\sum_{j=1}^N \Vert \by_k-A\nabla v\Vert^2_{A^{-1},\om_k} \right)^{1/2}\Vert \nabla w\Vert_{A,\Omega} \\
&+ \left( \displaystyle\sum_{j=1}^N C_{\rm P, \max}^2 \Vert \dvg \by_k+f \Vert^2_{\om_k} \right)^{1/2} \Vert \nabla w\Vert_{\Omega}\\
& +  \displaystyle\sqrt{E_{\max}} \left(\displaystyle\sum_{\gamma_{kj}\in\mathcal{E}_1}
\Vert (\by_k-\by_j)\cdot \bn_{kj}\Vert_{\gamma_{kj}}^2\beta_{kj}^2 \right)^{1/2} \Vert \nabla w\Vert_{\Omega}.
\end{array}
$$
Finally, setting $w=u-v$ and noting that \eqref{pd_A}
implies
\begin{equation}\label{eq:D10}
\Vert \nabla w\Vert^2 \le \frac{1}{C_{\min}}\int_{\Omega} A\nabla w\cdot \nabla w dx
\end{equation}
we arrive at the  a posteriori error estimate naturally
adapted to approximations   computed via the Schwarz alternating method according to Algorithm~\ref{alg1}.

\begin{theorem}\label{thm:main_result}
Let $v \in H^1(\Omega)$ be an approximation of the exact solution to the boundary-value
problem~\eqref{eq:1a}--\eqref{eq:1b} that  satisfies the boundary condition and let the flus field $\by=(\by_1,\ldots,\by_N) \in H(\bigcup_j \om_j,\dvg)$ satisfy the
constraints A.\ref{ass1} and~A.\ref{ass2}. Then the following estimate holds:
\begin{equation}\label{eq:D11}
\begin{array}{rl}
 \displaystyle\Vert \nabla(u-v)\Vert_{A,\Omega}& 
  \le
  \left( \displaystyle\sum_{k=1}^N \Vert \by_k-A\nabla v\Vert^2_{A^{-1},\om_k} \right)^{1/2} \\ &
  + \displaystyle \frac{1}{\sqrt{C_{\min}}}\left[ C_{\rm P,\max} \left( \displaystyle\sum_{k=1}^N \Vert \dvg \by_k+f \Vert^2_{\om_k} \right)^{1/2} \right. \\ &
  \left. + \displaystyle\sqrt{E_{\max}} \left(\displaystyle\sum_{\gamma_{kj}\in\mathcal{E}_1}\beta_{kj}^2 \Vert (\by_k-\by_j)\cdot \bn_{kj}\Vert_{\gamma_{kj}}^2 \right)^{1/2} \right].
\end{array}
\end{equation}
\end{theorem}

\begin{remark}\label{remark:estimate}
It is convenient to represent the estimate in a somewhat different
form. We square both parts of
\eqref{eq:D11}, apply Young's inequality with positive
$\varepsilon_1$, $\varepsilon_2$, and $\varepsilon_3$,
and obtain the upper bound
\begin{equation}\label{eq:D12}
\begin{array}{rl}
 \displaystyle\Vert \nabla(u-v)\Vert_{A,\Omega}^2& 
  \le
  \alpha_1 \displaystyle\sum_{k=1}^N \Vert \by_k-A\nabla v\Vert^2_{A^{-1},\om_k} \\ &
  + \alpha_2 \displaystyle\sum_{k=1}^N \Vert \dvg \by_k+f \Vert^2_{\om_k} \\ &
  + \alpha_3 \displaystyle\sum_{\gamma_{kj}\in
  \mathcal{E}_1}\beta_{kj}^2 \Vert (\by_k-\by_j)\cdot \bn_{kj}\Vert_{\gamma_{kj}}^2 
  =:M_{\oplus}(\by,v,f;\varepsilon_1,\varepsilon_2,\varepsilon_3),
\end{array}
\end{equation}
where
\be
&&\alpha_1:=1+\varepsilon_1+\varepsilon_2,\quad
\alpha_2:=(1+\varepsilon_1^{-1}+\varepsilon_3) C_{\rm P,\max}^2/C_{\min},\\
&&
\alpha_3 :=(1+\varepsilon_2^{-1}+\varepsilon_3^{-1}) E_{\max}/C_{\min}.
\ee
The right hand side of (\ref{eq:D12}) is composed of quadratic functionals which is more
convenient in the process of determining suitable flux $\by:=(\by_1,\by_2,\ldots,\by_N)$.
The constants $\varepsilon_i >0$ for $i=1,2,3$ in the definition of $\alpha_i$, $i=1,2,3$,
can be chosen independently from each other and such that they minimize the majorant
$M_{\oplus}(\by,v,f;\varepsilon_1,\varepsilon_2,\varepsilon_3)$.
\end{remark}

%%%%%%%%%%%%%%%%%%%%
\subsection{Computation of admissible flux fields}\label{subsec:admissible_fluxes}

In this subsection, we  address the problem of computing an admissible flux field $\by$ from a flux field
$\tilde{\by} \in H(\bigcup_j \om_j,\dvg)$, which may violate the conditions A.\ref{ass1} and~A.\ref{ass2}.

In the following we outline the procedure for a piecewise linear $H^1$ vector field $v$ obtained from
using $P_1$ finite element approximations in the subdomain solves of the domain decomposition method.
To begin with, we define
\begin{equation}\label{eq:D13}
\tilde{\by}_k:= G_k (A \nabla v) \vert_{\omega_k}, \quad k=1,2,\ldots,N,
\end{equation}
where $G_k$ is the gradient averaging operator acting on the gradient of a conforming, piecewise (elementwise)
linear approximation of the solution on the basic subdomain $\omega_k$. The resulting piecewise linear vector field
$\tilde{\by}_k$ satisfies $\tilde{\by}_k \in H^1(\omega_k) \subset H(\om_k,\dvg)$. As a consequence, the vector
field $\tilde{\by}:=(\tilde{\by}_1,\tilde{\by}_2,\ldots,\tilde{\by}_N)$ is a broken $H(\dvg)$ field,
i.e., $\tilde{\by} \in H(\bigcup_j \om_j,\dvg)$, and, may have discontinuous normal component on the
interfaces $\gamma_{kj}\in\mathcal{E}_1$. Moreover, in general, the field $\tilde{\by}$ will not satify the
Assumptions~A.\ref{ass1} and~A.\ref{ass2}, which are required to apply our theory.  

In order to overcome this difficulty,
we construct a new vector field $\by=(\by_1,\ldots,\by_N)$ from $\tilde{\by}$ by adding on each subdomain~$\omega_k$
a proper corrector $\bq_k \in H(\om_k,\dvg)$ to $\tilde{\by}_k$, i.e., we define
\begin{equation}\label{eq:D14}
\by_k := \tilde{\by}_k + \bq_k, \quad k=1,2,\ldots,N,
\end{equation}
and herewith the broken $H(\dvg)$ field $\by := \tilde{\by} + \bq$, i.e., $\by \in H(\bigcup_j \om_j,\dvg)$,
the latter satisfying Assumptions~A.\ref{ass1} and~A.\ref{ass2}.
In fact, the conditions to be satisfied are as follows:
\ben
\label{c1}
&&\int_{\om_k}\dvg(\widetilde\by_k+\bq_k) dx=0,\quad k,j=1,2,...,N,\\
\label{c2}
&&\int_{\gamma_{kj}}(\tilde\by_k-\tilde\by_j+\bq_k-\bq_j)\cdot \bn_{kj} d\gamma=0.
\een
We consider the simplest class of correctors, which are fully defined by constant normal fluxes
on the edges of a polygonal cell $\omega_k$. Extensions of the vector field $\bq$ inside $\omega_k$
can be obtained by Raviart-Thomas elements of the lowest order. If $\om_k$ is a simplex, then the corresponding
extension corresponds to the RT$^0$ element. If $\om_k$ is a polygonal cell with $t$ edges, then we introduce 
$t-1$ nonintersecting diagonals and use RT$^0$ elements in the emerging simplexes.

Overall, the correction field $\bq$ can be viewed as an approximation generated on a coarse
mesh ${\mathcal T}_\om$ associated with the subdomains $\om_k$. Let ${\bf Q}_N$ denote
the corresponding finite dimensional space. 
It should be outlined that the dimensionality of ${\bf Q}_N$ is small and fixed (it does not depend
on the subspaces used by the DDM in the subdomains $\Omega_i$). For example, if $d=2$ and
all $\omega_i$, $i=1,2,...,N$ are triangles, then ${\rm dim} {\bf Q}_N=3N$.

First of all we need to verify that {\em ${\rm dim} {\bf Q}_N$ is large enough for satisfying the
conditions}~(\ref{c1})--(\ref{c2}). Let $N_v$ and $N_f$ denote the number of vertices and faces
in $\mathcal T_\om$, respectively.
If $N_{fD}$ denotes the number of faces on the Dirichlet boundary, then we need to satisfy $N$
conditions of the form~(\ref{c1}) and $N_f-N_{fD}$ conditions~(\ref{c2}).
Hence the general condition is
\ben
\label{cond1}
{\rm dim} {\bf Q}_N+N_{fD}\geq \,N+N_f.
\een
Assume that all the cells $\omega_i$ have $\ell$ faces. Then ${\rm dim} {\bf Q}_N=N\ell$
and (\ref{cond1}) reads
\ben
\label{cond2}
N(\ell-1)+N_{fD}\geq\,N_f
\een
If $d=2$, then, from the Euler identity for planar graphs, it follows that
\be
N_v+(N+1)-N_f=2, \quad \mbox{or, equivalently,} \quad  N_f=N_v+N-1 .
\ee
Hence, we arrive at the simple compatibility condition
\ben
\label{cond3}
N(\ell-2)+N_{fD}\geq N_v-1.
\een

Let us illustrate (\ref{cond2})--(\ref{cond3}) with several examples.
\begin{figure}[h!]
\caption{Coarse mesh ${\mathcal T}_\om$.}
\label{figure:comp_cond}
\vspace{1ex}
\begin{tikzpicture}[scale=1.5] % picture (a)
\coordinate (1) at (0,0); \coordinate (2) at (1,0); \coordinate (3) at (2,0); 
\coordinate (5) at (0,1); \coordinate (6) at (1,1);\coordinate (7) at (2,1);
\coordinate (9) at (0,2); \coordinate (10) at (1,2); 
\node at (1){$\bullet$};\node at (2){$\bullet$};
\node at (3){$\bullet$};%\node at (4){$\bullet$};
\node at (5){$\bullet$};\node at (6){$\bullet$};
\node at (7){$\bullet$};%\node at (8){$\bullet$};
\node at (9){$\bullet$};\node at (10){$\bullet$};
%\node[below] at (1,4){$G$};
\node at (0.4,1.68){$\omega_1$};
\node at (0.6,1.3){$\omega_2$};
\node at (0.4,0.68){$\omega_3$};
\node at (0.6,0.25){$\omega_4$};
\node at (1.4,0.68){$\omega_5$};
\node at (1.6,0.25){$\omega_6$};
\node at (1,-0.5){(a)};
%\node at (1.5,0.5){$\omega_8$};
%\node at (2.5,0.5){$\omega_9$};
\draw  (1)--(2)--(3); % draw lines
\draw (5)--(6)--(7);
\draw  (9)--(10);
\draw (1)--(5)--(9); % draw lines
\draw (2)--(6)--(10);
\draw (3)--(7);\draw (1)--(6);
\draw (2)--(7);\draw (5)--(10);
\end{tikzpicture}
\begin{tikzpicture}[scale=0.69] % picture (b)
\coordinate (1) at (0,0); \coordinate (2) at (1,0); %declare a point and its coordinates
\coordinate (3) at (2,0); \coordinate (4) at (3,0);
\coordinate (5) at (0,1); \coordinate (6) at (1,1); %declare a point and its coordinates
\coordinate (7) at (2,1); \coordinate (8) at (3,1);
\coordinate (9) at (0,2); \coordinate (10) at (1,2); %declare a point and its coordinates
\coordinate (11) at (2,2); \coordinate (12) at (3,2);
\coordinate (13) at (0,3); \coordinate (14) at (1,3); %declare a point and its coordinates
\coordinate (15) at (2,3); \coordinate (16) at (3,3);
\coordinate (17) at (4,0); \coordinate (18) at (4,1); %declare a point and its coordinates
\coordinate (19) at (4,2); \coordinate (20) at (4,3);
\node at (0,0){$\bullet$}; % put node marked by a bold dot
%\node at (0,0) {$\circ$}; % put node marked by a circle
%\node at (0,0) {$\Box$}; % put node marked by box
%\node at (0,0) {$\star$}; % put node marked by star
\node at (1){$\bullet$};\node at (2){$\bullet$};
\node at (3){$\bullet$};\node at (4){$\bullet$};
\node at (5){$\bullet$};\node at (6){$\bullet$};
\node at (7){$\bullet$};\node at (8){$\bullet$};
\node at (9){$\bullet$};\node at (10){$\bullet$};
\node at (11){$\bullet$};\node at (12){$\bullet$};
\node at (13){$\bullet$};\node at (14){$\bullet$};
\node at (15){$\bullet$};\node at (16){$\bullet$};
\node at (17){$\bullet$};\node at (18){$\bullet$};
\node at (19){$\bullet$};\node at (20){$\bullet$};
\draw  (1)--(2)--(3)--(4)--(17); % draw lines
\draw (5)--(6)--(7)--(8)--(18);
\draw  (9)--(10)--(11)--(12)--(19);
%\draw  [ultra thick,red] (9)--(10)--(11)--(12);
\draw   (13)--(14)--(15)--(16)--(20);
\draw (1)--(5)--(9)--(13); % draw lines
\draw (2)--(6)--(10)--(14);
\draw (3)--(7)--(11)--(15);
\draw (4)--(8)--(12)--(16);
\draw (17)--(18)--(19)--(20);
\draw (1)--(6)--(11)--(16);
\draw (2)--(7)--(12)--(20);
\draw (5)--(10)--(15);\draw (9)--(14);
\draw (4)--(18);\draw (3)--(8)--(19);
\node at (-0.5,3){$n$};
\node at (4,-0.5){$m$};
\node at (2,-1.5){(b)};
\end{tikzpicture}
\begin{tikzpicture}[scale=0.69] % picture (c)
\coordinate (1) at (0,0); \coordinate (2) at (1,-0.5); %declare a point and its coordinates
\coordinate (3) at (2,-0.5); \coordinate (4) at (3,-0.5);
\coordinate (5) at (0.3,1); \coordinate (6) at (1,0.5); %declare a point and its coordinates
\coordinate (7) at (2,0.5); \coordinate (8) at (3,0.5);
\coordinate (9) at (0.3,2); \coordinate (10) at (1,2); %declare a point and its coordinates
\coordinate (11) at (2,2); \coordinate (12) at (3,2);
\coordinate (13) at (0,3); \coordinate (14) at (1,3.5); %declare a point and its coordinates
\coordinate (15) at (2,3.5); \coordinate (16) at (3,3.5);
\coordinate (17) at (4,0); \coordinate (18) at (4,1); %declare a point and its coordinates
\coordinate (19) at (4,2); \coordinate (20) at (4,3);
\node at (0,0){$\bullet$}; 
\node at (1){$\bullet$};\node at (2){$\bullet$};
\node at (3){$\bullet$};\node at (4){$\bullet$};
\node at (5){$\bullet$};\node at (6){$\bullet$};
\node at (7){$\bullet$};\node at (8){$\bullet$};
\node at (9){$\bullet$};\node at (10){$\bullet$};
\node at (11){$\bullet$};\node at (12){$\bullet$};
\node at (13){$\bullet$};\node at (14){$\bullet$};
\node at (15){$\bullet$};\node at (16){$\bullet$};
\node at (17){$\bullet$};\node at (18){$\bullet$};
\node at (19){$\bullet$};\node at (20){$\bullet$};
\node at (1.5,1.5){$\omega_k$};
\draw  (1)--(2)--(3)--(4)--(17); % draw lines
\draw (5)--(6)--(7)--(8)--(18);
\draw  (9)--(10)--(11)--(12)--(19);
%\draw  [ultra thick,red] (9)--(10)--(11)--(12);
\draw   (13)--(14)--(15)--(16)--(20);
\draw (1)--(5)--(9)--(13); % draw lines
\draw (2)--(6)--(10)--(14);
\draw (3)--(7)--(11)--(15);
\draw (4)--(8)--(12)--(16);
\draw (17)--(18)--(19)--(20);
\node at (-0.5,3){$n$};
\node at (4,-0.5){$m$};
\node at (2,-1.5){(c)};
\end{tikzpicture}
\end{figure}
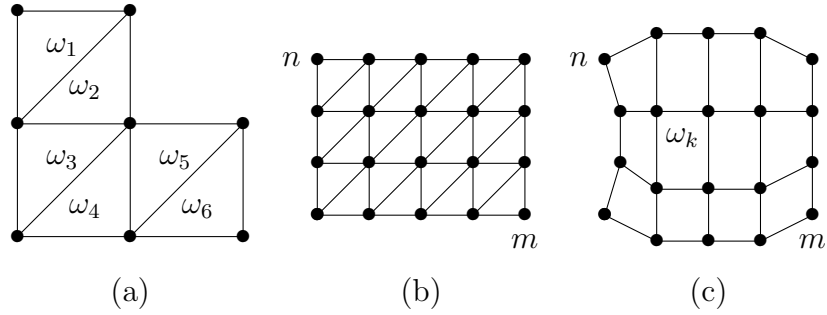
It is easy to see that for the mesh $\mathcal T_\om$ depicted in
Figure~\ref{figure:comp_cond}(a), $\ell=3$, $N=6$, and $N_v=8$,
so that (\ref{cond3}) holds only if $N_{fD}\geq 1$.
For a regular mesh depicted in Figure~\ref{figure:comp_cond}(b),
we have $N=2mn$, $N_v=(n+1)(m+1)$ and (\ref{cond3}) reads
\ben
\label{cond4}
mn+N_{fD}\geq m+n.
\een
If $N_{fD}=0$, then (\ref{cond4}) shows that the compatibility condition holds
only if $n>1$ and $m>1$. If $\mathcal T_\om$ contains quadrilateral cells
(as in~Figure~\ref{figure:comp_cond}(c)), then (\ref{cond3}) is also reduced
to~(\ref{cond4}).
These examples show that  the condition~(\ref{cond1}) should be verified for
a particular decomposition in order to avoid some special cases, in which the
topology of ${\mathcal T}_\om$ must be changed.  

If $t=N(\ell-2)+N_{fD}+1-N_v>0$,
then (\ref{cond3}) is satisfied, and, moreover, we have $t$ free parameters that
can be used to minimize local contributions
\be
\alpha_1 \Vert \by_1-A\nabla v\Vert^2_{A^{-1},\om_s} +
  \alpha_2 \displaystyle \Vert \dvg \by_1+f \Vert^2_{\om_s},\quad s=1,2,...,t
\ee
associated with $\omega_k$ adjacent to the Dirichlet boundary. However, this way
may be not optimal. The best possible correction field in ${\bf Q}_N$ can be constructed
by solving a subsidiary minimization problem, which we describe next.

For given functions $\tilde{\by} \in H(\bigcup_j \om_j,\dvg)$ and $v$, $f$,
and parameters $\varepsilon_1$, $\varepsilon_2$, $\varepsilon_3$, we minimize
\eqref{eq:D12} (inserting the definition~\eqref{eq:D14}), we respect to $\bq$
under the constraints A.\ref{ass1} and~A.\ref{ass2}, which are also formulated in terms of $\bq$.
Using the method of Langrange multipliers results in the following variational problem for the
computation of the correction field $\bq$ and Lagrange multiplier $\boldsymbol{\lambda}$:
Find $\bq:=(\bq_1,\bq_2,\ldots,\bq_N) \in  H(\bigcup_j \om_j,\dvg)$ and
$\boldsymbol{\lambda} = ({\lambda_{1,2},\ldots,\lambda_{N-1,N},\lambda_{1},\ldots,\lambda_N})$ such that
\begin{equation}\label{eq:D15abstract}
\mathcal{A}((\bq;\boldsymbol{\lambda}),(\br;\boldsymbol{\mu})) = \mathcal{F}((\br;\boldsymbol{\mu}))
\end{equation}
for all $\br:=(\br_1,\br_2,\ldots,\br_N) \in  H(\bigcup_j \om_j,\dvg)$ and
$\boldsymbol{\mu} = ({\mu_{1,2},\ldots,\mu_{N-1,N},\mu_{1},\ldots,\mu_N})$,
\begin{eqnarray}
\qquad \mathcal{A}((\bq;\boldsymbol{\lambda}),(\br;\boldsymbol{\mu}))&\hspace{-2ex}=\hspace{-2ex}&
\alpha_1 \displaystyle\sum_{k=1}^N \int_{\omega_k} \bq_k \cdot \br_k dx
+ \alpha_2  \displaystyle\sum_{k=1}^N \int_{\omega_k} \dvg \bq_k \, \dvg \br_k dx \label{eq:D15a} \\
&&
+ \alpha_3  \displaystyle\sum_{\gamma_{kj}\in\mathcal{E}_1}\beta_{kj}^2  \int_{\gamma_{kj}}
[(\bq_k-\bq_j) \cdot \bn_{kj}] [(\br_k-\br_j) \cdot \bn_{kj}] ds \nonumber \\
&&
+ \displaystyle\sum_{k=1}^N \int_{\omega_k} (\dvg \bq_k) \mu_k + (\dvg \br_k) \lambda_k dx \nonumber \\ 
&&
+ \vspace{-1ex} \displaystyle\sum_{\gamma_{kj}\in\mathcal{E}_1} \int_{\gamma_{kj}} [(\bq_k-\bq_j) \cdot \bn_{kj}] \mu_{k,j}
+   [(\br_k-\br_j) \cdot \bn_{kj}]  \lambda_{k,j} ds , \nonumber
\end{eqnarray}
\begin{eqnarray}
\qquad \mathcal{F}((\br;\boldsymbol{\mu}))
&=& \alpha_1 \displaystyle\sum_{k=1}^N \int_{\omega_k} (A \nabla v - \tilde{\by}_k) \cdot \br_k dx
- \alpha_2  \displaystyle\sum_{k=1}^N \int_{\omega_k} (f + \dvg \tilde{\by}_k) \, \dvg \br_k dx \label{eq:D15b} \\ 
&&
- \alpha_3  \displaystyle\sum_{\gamma_{kj}\in\mathcal{E}_1}\beta_{kj}^2 \int_{\gamma_{kj}}
[(\tilde{\by}_k-\tilde{\by}_j) \cdot \bn_{kj}] [(\br_k-\br_j) \cdot \bn_{kj}] ds \nonumber \\ 
&&
- \displaystyle\sum_{k=1}^N \int_{\omega_k} (f + \dvg \tilde{\by}_k) \mu_k dx \nonumber \\ 
&&
- \displaystyle\sum_{\gamma_{kj}\in\mathcal{E}_1} \int_{\gamma_{kj}}[(\tilde{\by}_k-\tilde{\by}_j)
\cdot \bn_{kj}]  \mu_{k,j} ds. \nonumber
\end{eqnarray}
Note that there is one Lagrange multiplier $\lambda_k$ per basic subdomain $\omega_k$ and
one Lagrange multiplier $\lambda_{k,j}$ for each interface $\gamma_{kj}\in\mathcal{E}_1$ in this
systems; the corresponding test functions are $\mu_k$ and $\mu_{k,j}$.

Once the corrector $\bq$ and thus the admissible flux field $\by$ has been found we can minimize
$M_{\oplus}(\tilde{\by} + \bq,v,f;\varepsilon_1,\varepsilon_2,\varepsilon_3)= M_{\oplus}(\varepsilon_1,\varepsilon_2,\varepsilon_3)$
with respect to $\boldsymbol \varepsilon:=(\varepsilon_1,\varepsilon_2,\varepsilon_3)$. Alternating minimization with respect to
one of the variables $\bq$ and $\boldsymbol \varepsilon$, keeping all others fixed, typically converges fast and we can stop
the procedure once a satifying result has been achieved.

The next section provides numerical evidence of the functionality of the proposed new methodology.

\section{Numerical evidence}\label{sec:numerics} 

In this final section we present some numerical tests to give numerical evidence that the new error
majorant~\eqref{eq:D11}, or, alternatively~\eqref{eq:D12}, can be applied successfully in combination
with DDM. We consider problem~\eqref{eq:1} with $A=I$ and $u_g$ chosen such the exact solution on
the L-shaped domain $\Omega=\Omega_1 \cup \Omega_2 := ((0,1)\times(0,2)) \cup ((0,2)\times(0,1))$,
as depicted in Figure~\ref{figure1}, is given by
\begin{equation}\label{eq:uex}
u=\frac{1}{\pi^2}\left(\sin(\pi x) \sin(\pi y)+\frac{1}{2}(1-\cos(\pi x)) (1-\cos(\pi y))\right).
\end{equation}

\begin{figure}[h]
\centering 
\hspace{-0.5cm}\subfigure[approximation~$\overline{v}^2$]{\label{fig:v_2}
\includegraphics[width=70mm]{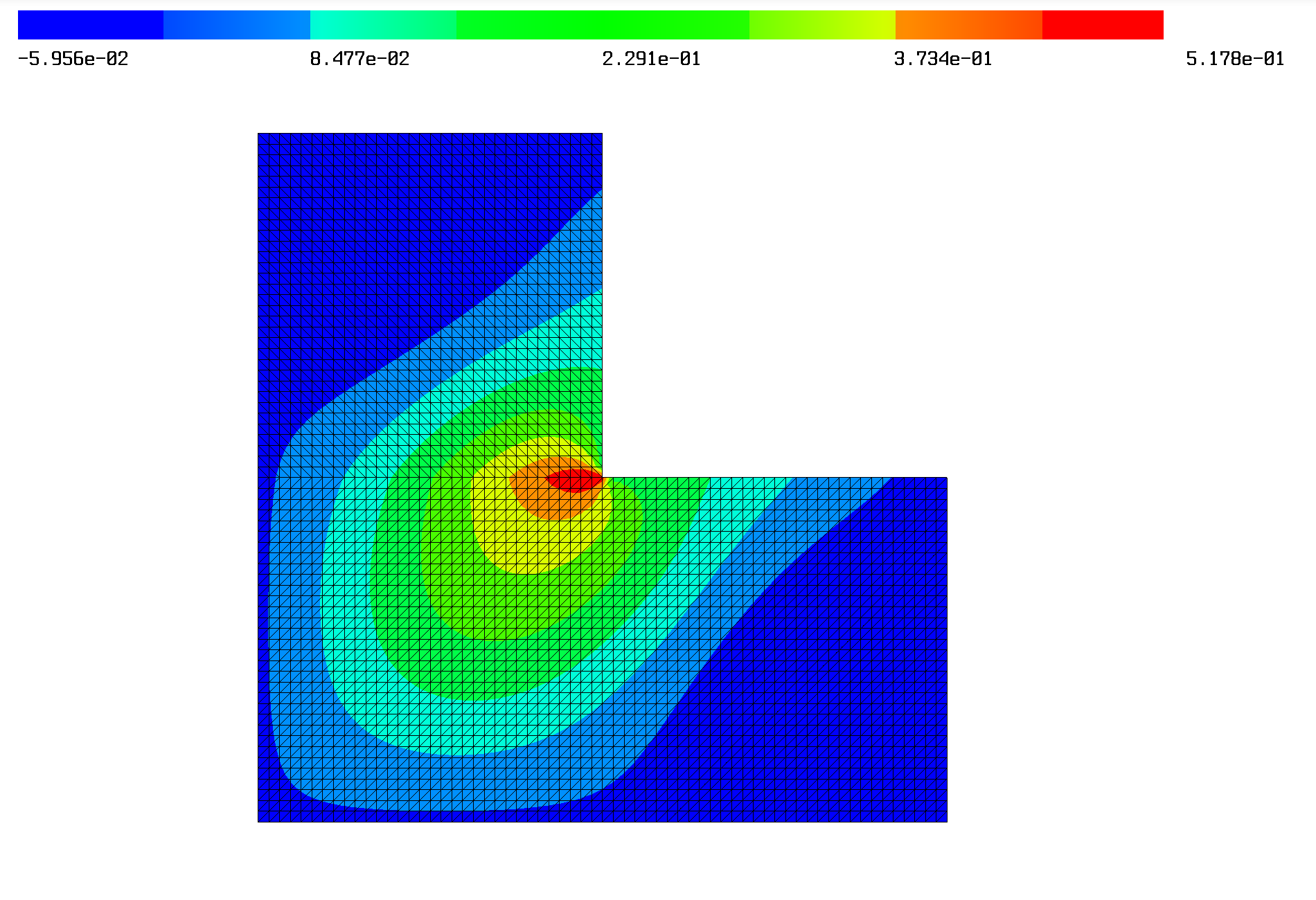}}\hspace{6ex}
\hspace{-0.5cm}\subfigure[approximation~$\overline{v}^4$]{\label{fig:v_4}
\includegraphics[width=70mm]{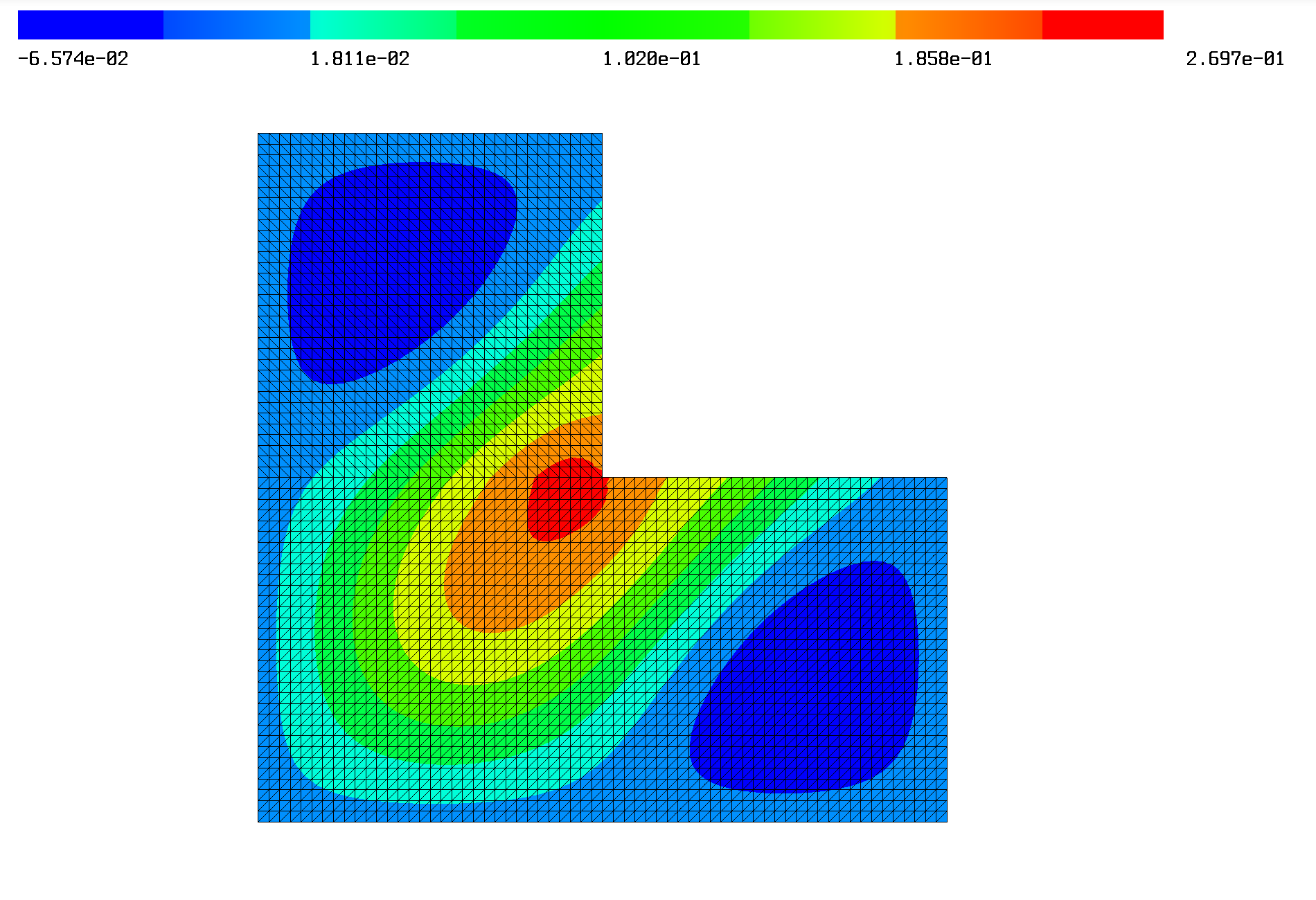}}\\
\hspace{-0.5cm}\subfigure[approximation~$\overline{v}^6$]{\label{fig:v_6}
\includegraphics[width=70mm]{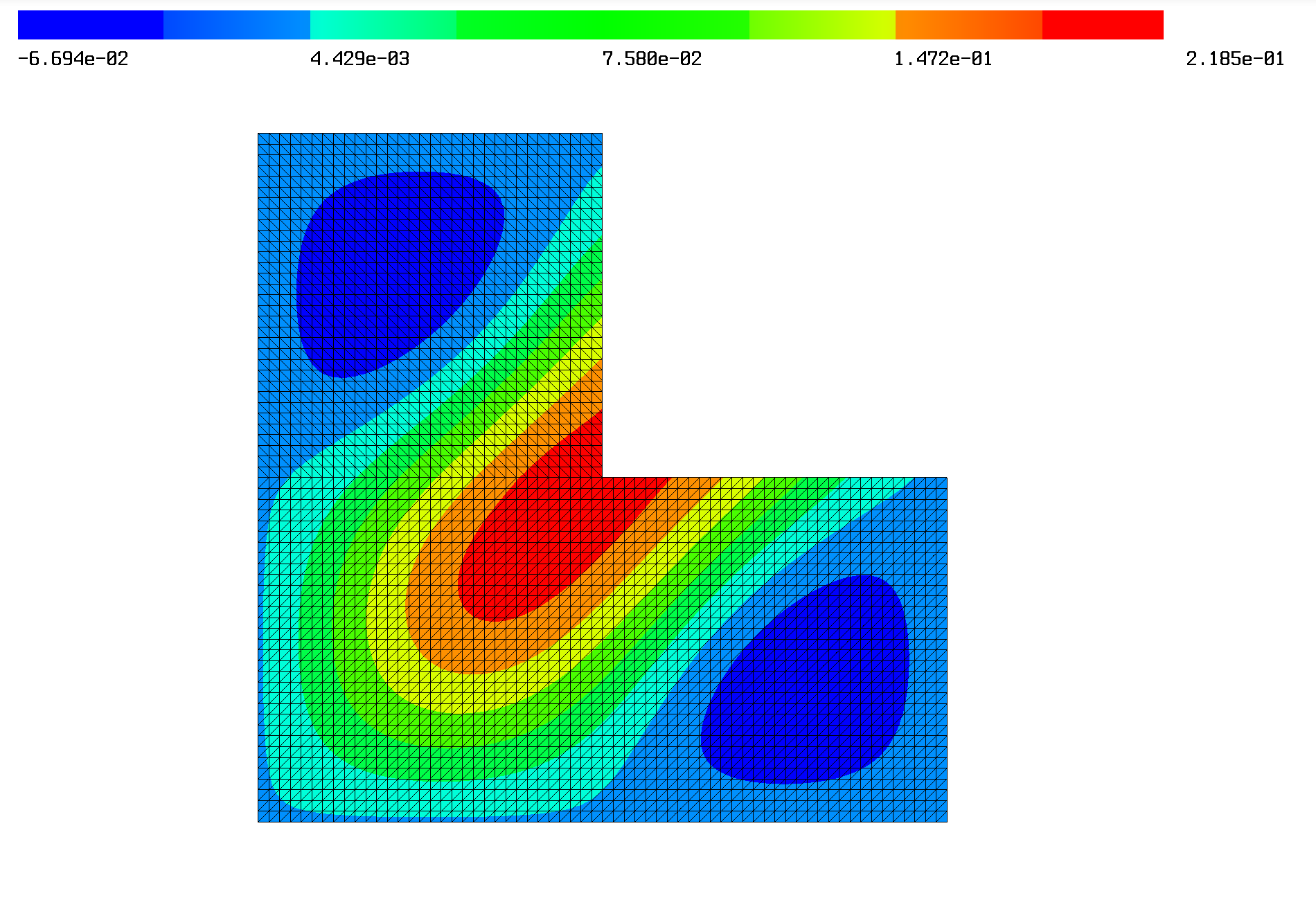}}\hspace{6ex}
\hspace{-0.5cm}\subfigure[approximation~$\overline{v}^8$]{\label{fig:v_8}
\includegraphics[width=70mm]{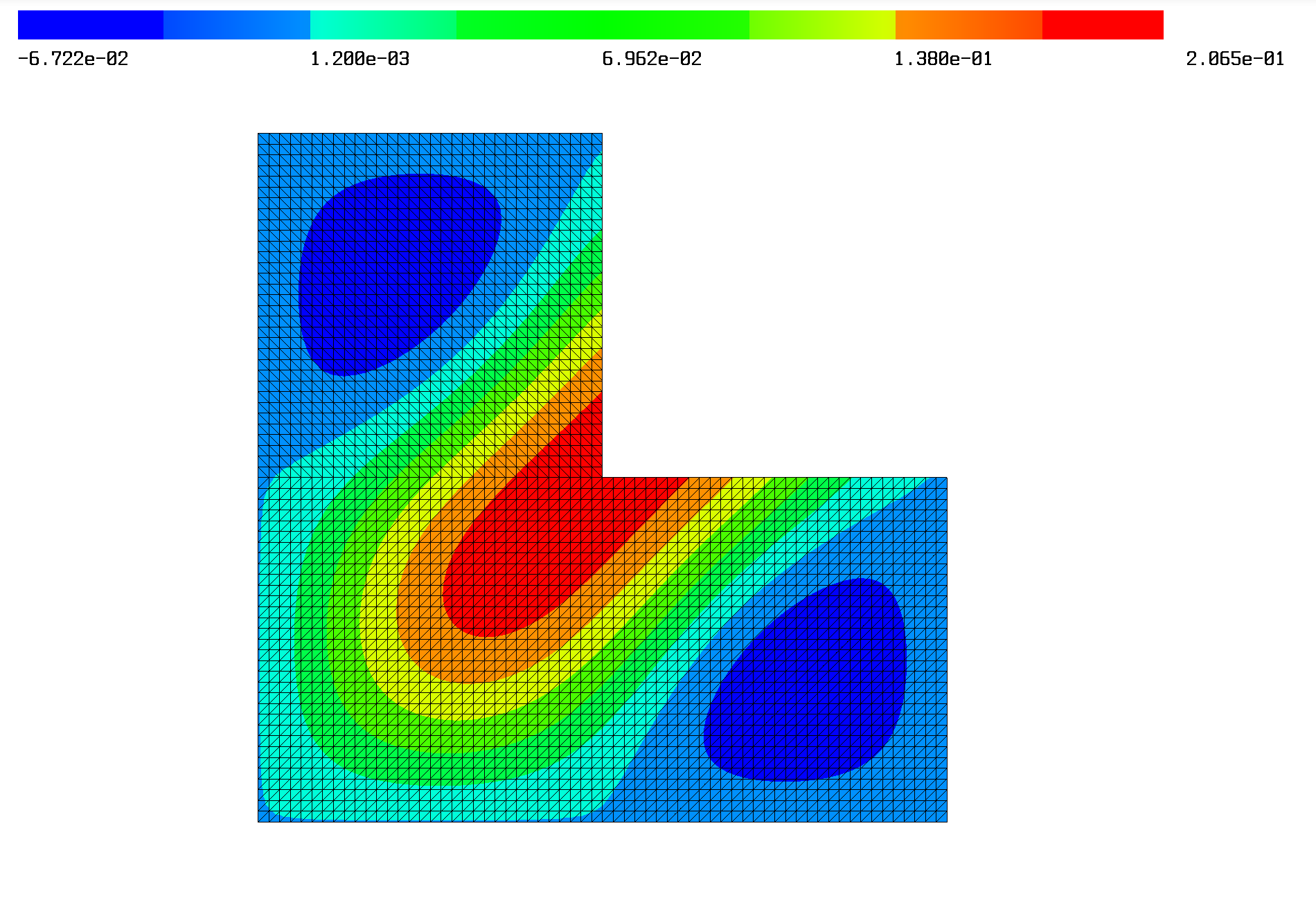}}\\
\caption{Global numerical approximation after iterations $2,4,6$
and $8$ of Algorithm~\ref{alg1}}\label{figure:approx_sol}.
\end{figure}

The global numerical approximations~$\overline{v}^n$, see~\eqref{glob_approx},
generated by Algorithm~\ref{alg1} in iterations $n=2,4,6,8$ are illustrated in
Figure~\ref{figure:approx_sol}, where $\overline{v}^n$ stems from using a standard
conforming finite element method (lowest-order Courant elements) to solve the subdomain
problems. All computations were performed using the NGSolve Finite Element Library
(http://sourceforge.net/projects/ngsolve).

We evaluated the contributions of different parts of the majorant~\eqref{eq:D12},
which we denote by $M_1^2$, $M_2^2$, and $M_3^2$, i.e., splitting the majorant according to
$$
M_{\oplus}^2:=M_1^2+M_2^2+M_3^2,
$$
where
\begin{subequations}\label{eq:parts}
\begin{eqnarray}
M_1^2 &:=&\alpha_1 \displaystyle\sum_{k=1}^N \Vert \by_k-A\nabla v\Vert^2_{A^{-1},\om_k}, \label{part1} \\
M_2^2 &:=& \alpha_2 \displaystyle\sum_{k=1}^N \Vert \dvg \by_k+f \Vert^2_{\om_k},  \label{part2} \\
M_3^2 &:=& \alpha_3 \displaystyle\sum_{\gamma_{kj}\in
  \mathcal{E}_1}\beta_{kj}^2 \Vert (\by_k-\by_j)\cdot \bn_{kj}\Vert_{\gamma_{kj}}^2.  \label{part3}
\end{eqnarray}
\end{subequations}
In the computations we did not perform any optimization of the parameters $\varepsilon_i$, $i=1,2,3$,
but used $\varepsilon_1=\varepsilon_2=\varepsilon_3=1$ instead. Further used in the computations
are the constants $C_{\rm P,\max}=\sqrt{2}/\pi, C_{\min}=1, E_{\max}=2$, $\alpha_i$, $i=1,2,3$, as in
Remark~\ref{remark:estimate}, and $\beta_{12}=\beta_{23}=1/\sqrt{\pi \tanh(\pi)}$.

Tables~\ref{tab1}--\ref{tab2} are for different mesh size. The quantities summarized in Table~\ref{tab1}
were computed with approximate solutions $v=\overline{v}^n$ and a flux corrector $\bq$ on the same mesh
with mesh size $h=H$ whereas the results in Table~\ref{tab2} were obtained with a flux corrector on a
coarser mesh with mesh size $H>h$.

As we can see, the efficiency index
$$
I_{\rm eff} := \frac{M_{\oplus}(\by,v,f;1,1,1)}{\displaystyle\Vert \nabla(u-v)\Vert_{A,\Omega}}
$$
in the former case is very good even without optimization of the parameters $\varepsilon_i$, $i=1,2,3$,
and gets much worse in the latter case. This suggests the strategy to improve a flux corrector that has
been computed from a global problem on a coarse mesh with mesh size $H$ by solving local subdomain
problems on a fine mesh with mesh size~$h$. 
\begin{table}[h!]
\caption{Contributions to error majorant~\eqref{eq:D12}; Evaluation after 16 iterations of Algorithm~\ref{alg1}; 
Corrector computed on fine mesh, i.e., $H=h$.}
\label{tab1}
\begin{tabular}{l|ccc|c|c}
$h$ & $M_1^2$ & $M_2^2$ & $M_3^2$ & $M_{\oplus}^2$ &  $I_{\rm eff}$ \\ \hline
1/4   & 8.28e-2 & 1.49e-1 & 5.45e-4 & 2.32e-1 & 2.96 \\
1/8   & 2.23e-2 & 3.80e-2 & 2.16e-4 & 6.05e-2 & 2.96 \\
1/16 & 5.91e-3 & 9.56e-3 & 6.69e-5 & 1.55e-2 & 2.98 \\
1/32 & 1.62e-3 & 2.39e-3 & 1.82e-5 & 4.03e-3 & 3.04 \\
1/64 & 4.80e-4 & 5.98e-4 & 4.69e-6 & 1.08e-3 & 3.15 \\
\end{tabular}
\vspace{2ex}
\end{table}
\begin{table}[h!]
\caption{Contributions to error majorant~\eqref{eq:D12}; Evaluation after 16 iterations of Algorithm~\ref{alg1}; 
Approximations $v=\overline{v}^{n}$ computed on fine mesh with mesh size $h=1/64$, corrector $\bq$ on
coarser mesh with mesh size $H>h$.}
\label{tab2}
\begin{tabular}{l|ccc|c|c}
$H$ & $M_1^2$ & $M_2^2$ & $M_3^2$ & $M_{\oplus}^2$ &  $I_{\rm eff}$ \\ \hline
1/4   & 2.75e-3 & 4.62e-1 & 1.21e-3 & 4.66e-1 & 6.53e1 \\
1/8   & 1.74e-3 & 4.06e-1 & 3.07e-4 & 4.08e-1 & 6.11e1 \\
1/16 & 1.09e-3 & 3.17e-1 & 7.61e-5 & 3.18e-1 & 5.40e1 \\
1/32 & 7.17e-4 & 1.95e-1 & 1.88e-5 & 1.96e-1 & 4.23e1 \\
\end{tabular}
\vspace{2ex}
\end{table}

Next, we tested how the parts $M_1^2$, $M_2^2$, and $M_3^2$ of the majorant $M_{\oplus}^2$
change throughout the Schwarz alternating iterative prozess, i.e., for increasing iteration count $n$.
The results are presented in Table~\ref{tab3}.

\begin{table}[h!]
\caption{Contributions to error majorant~\eqref{eq:D12}; Increasing number of iterations of Algorithm~\ref{alg1}; 
Corrector computed on fine mesh, i.e., $H=h$.}
\label{tab3}
\begin{tabular}{l|ccc|c|c}
$n$ & $M_1^2$ & $M_2^2$ & $M_3^2$ & $M_{\oplus}^2$ &  $I_{\rm eff}$ \\ \hline
2 & 1.26e0  & 4.27e-2 & 1.56e-1 & 1.46e0  & 1.70 \\
4 & 5.91e-2 & 1.99e-3 & 8.63e-3 & 6.97e-2 & 1.72 \\
6 & 3.57e-3 & 6.59e-4 & 4.65e-4 & 4.70e-3 & 1.86 \\
8 & 6.55e-4 & 6.01e-4 & 2.57e-5 & 1.28e-3 & 2.66 \\
\end{tabular}
\vspace{2ex}
\end{table}

Finally, we tested how the parts $M_1^2$ and $M_2^2$, which can be evaluated on different basic
subdomains independently from each other, change as the number of iterations~$n$ increases.
As one can observe, see Table~\ref{tab4}, especially $M_1$ indicates well which basic subdomains
should be processed next by the Schwarz method in order to reduce the global error majorant effectively.
\begin{table}[h!]
\caption{Contributions to volume terms of the error majorant~\eqref{eq:D12} from different basic subdomains;
Increasing number of iterations of Algorithm~\ref{alg1}; Corrector computed on fine mesh, i.e., $H=h$.}
\label{tab4}
\begin{tabular}{l|ccc|ccc}
$n$ & $M_{1,\omega_1}^2$ & $M_{1,\omega_2}^2$ & $M_{1,\omega_3}^2$ & 
	  $M_{2,\omega_1}^2$ & $M_{2,\omega_2}^2$ & $M_{2,\omega_3}^2$\\ \hline
2 & 7.94e-1 & 3.10e-1 & 1.60e-1 & 1.80e-2 & 2.22e-2 & 2.50e-3 \\
3 & 3.59e-2 & 6.18e-2 & 1.71e-1 & 7.60e-4 & 3.99e-3 & 3.10e-3 \\
4 & 3.76e-2 & 1.31e-2 & 8.42e-3 & 8.78e-4 & 8.89e-4 & 2.23e-4 \\
7 & 3.01e-4 & 3.67e-4 & 5.54e-4 & 3.14e-4 & 1.58e-4 & 1.39e-4 \\
8 & 3.03e-4 & 2.55e-4 & 9.66e-5 & 3.15e-4 & 1.52e-4 & 1.34e-4 \\
\end{tabular}
\vspace{2ex}
\end{table}
These experiments confirm that the new localized a posteriori error majorant provided
by Theorem~\ref{thm:main_result} (and estimate~\eqref{eq:D12}) has great potential
to become a powerful tool for reliable and cost-efficient iterative solution methods for
(elliptic) PDEs by domain decomposition methods. 

Future investigations will deal with refining the proposed approach and adapting it to
various classes of problems.

{\bf Acknowledgements.}
The first author would like to thank Philip Lederer for his support with the high performance
multiphysics finite element software Netgen/NGSolve that served as a platform to conduct
the numerical tests.

%\bibliographystyle{plain}
%\bibliography{references}

\begin{thebibliography}{10}

\bibitem{BrSh}
D.~Braess and J.~Sch\"{o}berl.
\newblock Equilibrated residual error estimator for edge elements.
\newblock {\em Math. Comp.}, 77(262):651--672, 2008.

\bibitem{Bramble_etal1991}
J.H.~Bramble, J.E.~Pasciak, J.P.~Wang, and J.~Xu.
\newblock Convergence estimates for product iterative methods with applications
  to domain decomposition.
\newblock {\em Math. Comp.}, 57(195):1--21, 1991.

\bibitem{Do03}
C.R.~Dohrmann.
\newblock {A preconditioner for substructuring based on constrained energy
  minimization}.
\newblock {\em SIAM J. Sci. Comput.}, 25(1):246--258, 2003.

\bibitem{DoleanJolivetNataf2015}
V.~Dolean, P.~Jolivet, and F.~Nataf.
\newblock {\em An introduction to domain decomposition methods}.
\newblock Society for Industrial and Applied Mathematics (SIAM), Philadelphia,
  PA, 2015.
\newblock Algorithms, theory, and parallel implementation.

\bibitem{EGLW12}
Y.~Efendiev, J.~Galvis, R.~Lazarov, and J.~Willems.
\newblock {Robust domain decomposition preconditioners for abstract symmetric
  positive definite bilinear forms}.
\newblock {\em ESAIM: Mathematical Modelling and Numerical Analysis},
  46(05):1175--1199, 2012.

\bibitem{FLLPR01}
C.~Farhat, M.~Lesoinne, P.~LeTallec, K.~Pierson, and D.~Rixen.
\newblock {FETI-DP: A dual-primal unified FETI method. I. A faster alternative
  to the two-level FETI method}.
\newblock {\em Int. J. Numer. Methods Engrg.}, 50(7):1523--1544, 2001.

\bibitem{GE10a}
J.~Galvis and Y.~Efendiev.
\newblock {Domain decomposition preconditioners for multiscale flows in
  high-contrast media}.
\newblock {\em Multiscale Model. Simul.}, 8(4):1461--1483, 2010.

\bibitem{Ha03}
W.~Hackbusch.
\newblock {\em {Multi-Grid Methods and Applications}}.
\newblock Springer, Berlin Heidelberg, 2003.

\bibitem{KK}
L.V.~Kantorovich and V.I.~Krylov.
\newblock {\em Approximate methods of higher analysis}.
\newblock Interscience, 1964.

\bibitem{Ke}
D.W.~Kelly.
\newblock The self-equilibration of residuals and complementary a posteriori
  error estimates in the finite element method.
\newblock {\em Internat. J. Numer. Methods Engrg.}, 20(8):1491--1506, 1984.

\bibitem{Kr06}
J.~Kraus.
\newblock {Algebraic multilevel preconditioning of finite element matrices
  using local {S}chur complements}.
\newblock {\em Numer. Linear Algebra Appl.}, 13:49--70, 2006.

\bibitem{Kr12}
J.~Kraus.
\newblock {Additive Schur complement approximation and application to
  multilevel preconditioning}.
\newblock {\em SIAM J. Sci. Comput.}, 34:A2872--A2895, 2012.

\bibitem{KLM15}
J.~Kraus, M.~Lymbery, and S.~Margenov.
\newblock Auxiliary space multigrid method based on additive {S}chur complement
  approximation.
\newblock {\em Numer. Linear Algebra Appl.}, 22(6):965--986, 2015.

\bibitem{KLLMZ16}
J.~Kraus, R.~Lazarov, M.~Lymbery, S.~Margenov, and L.~Zikatanov.
\newblock Preconditioning heterogeneous {$H({\rm div})$} problems by additive
  {S}chur complement approximation and applications.
\newblock {\em SIAM J. Sci. Comput.}, 38(2):A875--A898, 2016.

\bibitem{K89}
Y.~Kuznetsov.
\newblock {Algebraic multigrid domain decomposition methods}.
\newblock {\em Sov. J. Numer. Anal. Math. Modelling.}, 4(5):351--379, 1989.

\bibitem{LaLe}
P.~Ladev\`eze and D.~Leguillon.
\newblock Error estimate procedure in the finite element method and
  applications.
\newblock {\em SIAM J. Numer. Anal.}, 20(3):485--509, 1983.

\bibitem{Lions_1978}
P.-L. Lions.
\newblock Interpr\'{e}tation stochastique de la m\'{e}thode altern\'{e}e de
  {S}chwarz.
\newblock {\em C. R. Acad. Sci. Paris}, 268:325--328, 1978.

\bibitem{Lions_1987}
P.-L. Lions.
\newblock On the {S}chwarz alternating method. {I}.
\newblock In R.~Glowinski, G.~H. Golub, G.~A. Meurant, and J.~P\'{e}riaux,
  editors, {\em First International Symposium on Domain Decomposition Methods
  for Partial Differential Equations}, pages 1--42, Paris, France, 1987. SIAM,
  Philadelphia, PA.

\bibitem{Lions_1989}
P.-L. Lions.
\newblock On the {S}chwarz alternating method. {II}.
\newblock In T.~Chan, R.~Glowinski, J.~P\'{e}riaux, and O.~Widlund, editors,
  {\em Second International Symposium on Domain Decomposition Methods for
  Partial Differential Equations}, pages 47--70, Paris, France, 1989. SIAM, Los
  Angeles, CA.

\bibitem{MaNeRe}
O.~Mali, P.~Neittaanm\"{a}ki, and S.~Repin.
\newblock {\em Accuracy verification methods}, volume~32 of {\em Computational
  Methods in Applied Sciences}.
\newblock Springer, Dordrecht, 2014.
\newblock Theory and algorithms.

\bibitem{Ma93}
J.~Mandel.
\newblock {Balancing domain decomposition}.
\newblock {\em Comm. Numer. Methods Engrg.}, 9(3):233--241, 1993.

\bibitem{MD03}
J.~Mandel and C.R.~Dohrmann.
\newblock {Convergence of a balancing domain decomposition by constraints and
  energy minimization}.
\newblock {\em Numer. Linear Algebra Appl.}, 10(7):639--659, 2003.

\bibitem{MDT05}
J.~Mandel, C.R.~Dohrmann, and R.~Tezaur.
\newblock {An algebraic theory for primal and dual substructuring methods by
  constraints}.
\newblock {\em Appl. Numer. Math.}, 54(2):167--193, 2005.

\bibitem{MS07}
J.~Mandel and B.~Soused{\'{\i}}k.
\newblock {Adaptive selection of face coarse degrees of freedom in the {BDDC}
  and {FETI-DP} iterative substructuring methods}.
\newblock {\em Comput. Methods Appl. Mech. Engrg.}, 196(8):1389--1399, 2007.

\bibitem{M08}
T.P.A.~Mathew.
\newblock {\em {Domain Decomposition Methods for the Numerical Solution of
  Partial Differential Equations}}.
\newblock Springer, Berlin Heidelberg, 2008.

\bibitem{MatsokinNepomnyashchikh1981}
A.M.~Matsokin and S.V.~Nepomnyashchikh.
\newblock Convergence of the {S}chwarz alternation method over subdomains
  without overlap.
\newblock In {\em Approximation and interpolation methods ({N}ovosibirsk,
  1980)}, pages 85--97. Akad. Nauk SSSR Sibirsk. Otdel., Vychisl. Tsentr,
  Novosibirsk, 1981.
\newblock Translated in Soviet J. Numer. Anal. Math. Modelling {{\bf{4}}}
  (1989), no. 6, 479--485.

\bibitem{MatsokinNepomnyashchikh1985}
A.M.~Matsokin and S.V.~Nepomnyashchikh.
\newblock The {S}chwarz alternation method in a subspace.
\newblock {\em Izv. Vyssh. Uchebn. Zaved. Mat.}, (10):61--66, 85, 1985.

\bibitem{Mikhlin}
S.G.~Mikhlin.
\newblock On the {S}chwarz algorithm.
\newblock {\em Doklady Akad. Nauk SSSR (N.S.)}, 77:569--571, 1951.

\bibitem{Mi}
S.G.~Mikhlin.
\newblock {\em Variational methods in mathematical physics}.
\newblock A Pergamon Press Book. The Macmillan Co., New York, 1964.
\newblock Translated by T. Boddington; editorial introduction by L. I. G.
  Chambers.

\bibitem{Nazarov_2015}
A.I~ Nazarov and S.~Repin.
\newblock Exact constants in {P}oincar\'{e} type inequalities for functions
  with zero mean boundary traces.
\newblock {\em Math Methods Appl Sci.}, 38:3195--3207, 2015.

\bibitem{Nev}
O.~Nevanlinna.
\newblock Remarks on {P}icard-{L}indel\"{o}f iteration. {II}.
\newblock {\em BIT}, 29(3):535--562, 1989.

\bibitem{Os}
A.~Ostrowski.
\newblock Les estimations des erreurs a posteriori dans les proc\'{e}d\'{e}s
  it\'{e}ratifs.
\newblock {\em C. R. Acad. Sci. Paris S\'{e}r. A-B}, 275:A275--A278, 1972.

\bibitem{ReGruyter}
S.~Repin.
\newblock A posteriori error estimation methods for partial differential
  equations.
\newblock In {\em Lectures on advanced computational methods in mechanics},
  volume~1 of {\em Radon Ser. Comput. Appl. Math.}, pages 161--226. Walter de
  Gruyter, Berlin, 2007.

\bibitem{Repin2000}
S.I.~Repin.
\newblock A posteriori error estimation for variational problems with uniformly
  convex functionals.
\newblock {\em Math. Comp.}, 69(230):481--500, 2000.

\bibitem{Sa94}
M.V.~Sarkis~Martins.
\newblock {\em Schwarz preconditioners for elliptic problems with discontinuous
  coefficients using conforming and non-conforming elements}.
\newblock ProQuest LLC, Ann Arbor, MI, 1994.
\newblock Thesis (Ph.D.)--New York University.

\bibitem{SVZ11}
R.~Scheichl, P.~Vassilevski, and L.~Zikatanov.
\newblock {Weak approximation properties of elliptic projections with
  functional constraints}.
\newblock {\em Multiscale Model. Simul.}, 9(4):1677--1699, 2011.

\bibitem{Schwartz}
H.A.~Schwarz.
\newblock Ueber einige {A}bbildungsaufgaben.
\newblock {\em J. Reine Angew. Math.}, 70:105--120, 1869.

\bibitem{Stoutemyer1973}
D.R.~Stoutemyer.
\newblock Numerical implementation of the {S}chwarz alternating procedure for
  elliptic partial differential equations.
\newblock {\em SIAM J. Numer. Anal.}, 10:308--326, 1973.

\bibitem{Sy}
J.L.~Synge.
\newblock The hypercircle method.
\newblock In {\em Studies in numerical analysis (papers in honour of
  {C}ornelius {L}anczos on the occasion of his 80th birthday)}, pages 201--217.
  1974.

\bibitem{dd_book}
A.~Toselli and O.~Widlund.
\newblock {\em Domain Decomposition Methods - Algorithms and Theory}.
\newblock Springer, Berlin Heidelberg, Germany, 2005.

\bibitem{TOS01}
U.~Trottenberg, C.W. Oosterlee, and A.~Sch{\"u}ller.
\newblock {\em {Multigrid}}.
\newblock Academic Press Inc., San Diego, CA, 2001.

\bibitem{V08}
P.~Vassilevski.
\newblock {\em {Multilevel Block Factorization Preconditioners: Matrix-based
  Analysis and Algorithms for Solving Finite Element Equations}}.
\newblock Springer, New York, 2008.

\bibitem{Ve3}
R.~Verf\"{u}rth.
\newblock {\em A review of a posteriori error estimation and adaptive mesh-refinement techniques}.
\newblock John Wiley and Sons Ltd, 1996.

\bibitem{Xu96}
J.~Xu.
\newblock {The auxiliary space method and optimal multigrid preconditioning
  techniques for unstructured grids}.
\newblock {\em Computing}, 56:215--235, 1996.

\bibitem{XZ02}
J.~Xu and L.~Zikatanov.
\newblock The method of alternating projections and the method of subspace
  corrections in {H}ilbert space.
\newblock {\em J. Amer. Math. Soc.}, 15(3):573--597, 2002.

\end{thebibliography}

\end{document}